\documentclass[11pt]{article}

\usepackage[latin1]{inputenc}\usepackage{epsfig}\usepackage{amsfonts}

\title{
  {\huge Fourier Theory on the Complex Plane V} \\
  Arbitrary-Parity Real Functions, \\
  Singular Generalized Functions and \\
  Locally Non-Integrable Functions }

\author{
  \Large Jorge L. deLyra \\
  Department of Mathematical Physics \\
  Physics Institute \\
  University of São Paulo }

\date{May 9, 2015}

\newcommand{\ii}{\mbox{\boldmath$\imath$}}

\newcommand{\e}[1]{\,{\rm e}^{#1}}
\newcommand{\blank}{\;\rule[-0.8ex]{0.8em}{0.05ex}\;}
\newcommand{\ldot}{\mbox{\Large$\cdot$}\!}
\newcommand{\Ev}{\mathfrak{E}}
\newcommand{\Od}{\mathfrak{O}}

\begin{document}\maketitle

\begin{abstract}
  \noindent
  A previously established correspondence between definite-parity real
  functions and inner analytic functions is generalized to real functions
  without definite parity properties. The set of inner analytic functions
  that corresponds to the set of all integrable real functions is then
  extended to include a set of singular ``generalized functions'' by the
  side of the integrable real functions. A general definition of these
  generalized functions is proposed and explored. The generalized
  functions are introduced loosely in the spirit of the Schwartz theory of
  distributions, and include the Dirac delta ``function'' and its
  derivatives of all orders. The inner analytic functions corresponding to
  this infinite set of singular real objects are given by means of a
  recursion relation. The set of inner analytic functions is then further
  extended to include a certain class of non-integrable real functions.
  The concept of integral-differential chains is used to help to integrate
  both the normal functions and the singular generalized functions
  seamlessly into a single structure. It does the same for the class of
  non-integrable real functions just mentioned. This extended set of
  generalized functions also includes arbitrary real linear combinations
  of all these real objects. An interesting connection with the Dirichlet
  problem on the unit disk is established and explored.
\end{abstract}

\section{Introduction}

In a previous paper~\cite{FTotCPI} we established a relation between
Definite-Parity (DP) real functions $f(\theta)$ defined on the interval
$[-\pi,\pi]$ and analytic functions $w(z)$ within the open unit disk of
the complex plane, as well as between the Fourier series of the DP real
functions and the complex power series of the corresponding analytic
functions, which we named ``inner analytic functions''. We refer the
reader to that paper for the detailed definition and discussion of many of
the concepts and notations we will use here. The DP real functions are
interpreted as the restrictions of these inner analytic functions to the
unit circle of the complex plane, and the Fourier series of the DP real
functions as the restrictions of the complex Taylor series of the
corresponding inner analytic functions to that same circle.

Within this context we showed, in a subsequent paper~\cite{FTotCPIV}, that
every DP real function within $[-\pi,\pi]$ that is absolutely integrable
is associated to a specific inner analytic function and can be recovered
from it almost everywhere in the limit from within the open unit disk to
the unit circle. This is true even if the corresponding Fourier series is
divergent, and the recovery of the DP real function can be executed using
only its sequence of Fourier coefficients $a_{k}$, which are thus seen to
uniquely characterize the DP real function almost everywhere. Let us
emphasize that, when we talk in this paper of the recovery of a real
function as the limit of an inner analytic function to the unit circle, we
always mean recovery almost everywhere, that is, except possibly in a
zero-measure subset of the domain.

On the other hand, if we consider the set of all possible inner analytic
functions, we realize that it includes much more that just those inner
analytic functions which are associated to the integrable DP real
functions. This is illustrated by the interesting and suggestive fact that
a radically singular object such as the Dirac delta ``function'' can also
be represented by an inner analytic function, as was shown in detail
in~\cite{FTotCPI}. As we will see, some non-integrable real functions can
also be represented by inner analytic functions. This at once poses the
question of what is the complete set of real objects on the unit circle
that corresponds to the set of all inner analytic functions within the
open unit disk.

Here we propose to define the objects within this set as {\em generalized
  functions}, loosely in the spirit of the Schwartz theory of
distributions~\cite{DistTheory}, which will include objects such as the
delta ``function'', its derivatives of arbitrarily high orders, and
possibly other singular objects. Although we will stop short of giving a
complete and detailed characterization of all possible such generalized
functions, we will show that many of the better-known ones are included in
the set. In fact, the inner analytic functions that correspond to the
delta ``function'' and to its derivatives of all orders will be exhibited
by a process of finite induction, leading to a simple algebraic recursion
relation. We will also show that many non-integrable functions are in the
set as well.

In order to better focus the analysis it may make good sense to impose
some limitations on the definition of the generalized functions. The most
general definition, as described above, and including both normal
functions and generalized functions under the heading of ``generalized
functions'', would be the following:

\vspace{3ex}

\noindent
\parbox{\textwidth} {\bf\boldmath The set of all generalized functions
  $f(\theta)$ with domain on the interval $[-\pi,\pi]$ is the set of the
  limits of the real and imaginary parts of the inner analytic functions
  $w(z)$ from within the open unit disk to the unit circle, whenever these
  limits exist at least almost everywhere.}

\vspace{3ex}

\noindent
In this definition all possible inner analytic functions, as they were
defined in~\cite{FTotCPI}, are included, and one may as well extend the
definition to all analytic functions within the open unit disk. The extra
conditions defining an inner analytic function $w(z)$, which were
introduced in~\cite{FTotCPI}, are that $w(0)=0$ and that $w(z)$ reduces to
a real function on the $(-1,1)$ interval of the real axis. These two extra
conditions were originally imposed mostly in order to simplify the
analysis. The first one is equivalent to the requirement that the real
functions be zero-average functions, which is a trivial limitation since
the addition of constant functions to the zero-average real functions,
thus lifting the limitation, is a trivial operation that has no
significant impact on any of the main results. The second one is
equivalent to the requirement that the real and imaginary parts of $w(z)$,
as well as the corresponding real functions, have definite parity
properties with respect to $\theta$. Since any real function $f(\theta)$
defined on $[-\pi,\pi]$ can be written in a unique way as a sum of an even
function and an odd function, we see that this requirement is not really a
limitation, and is used just to simplify the analysis.

The separation into even and odd parts can be easily applied to the
generalized functions as well as to normal functions. For example, the
Dirac delta ``function'' for an arbitrary singular point $\theta_{1}$ is
neither zero-average nor definite-parity, but the following zero-average
combination can be separated into a zero-average even part and an odd
part,

\noindent
\begin{eqnarray}\label{deltaparity}
  \delta(\theta-\theta_{1})-\frac{1}{2\pi}
  & = &
  \Ev\!\left[\delta(\theta-\theta_{1})-\frac{1}{2\pi}\right]
  +
  \Od\!\left[\delta(\theta-\theta_{1})-\frac{1}{2\pi}\right],
  \nonumber\\
  \Ev\!\left[\delta(\theta-\theta_{1})-\frac{1}{2\pi}\right]
  & = &
  \frac{\delta(\theta-\theta_{1})+\delta(\theta+\theta_{1})}{2}
  -
  \frac{1}{2\pi},
  \nonumber\\
  \Od\!\left[\delta(\theta-\theta_{1})-\frac{1}{2\pi}\right]
  & = &
  \frac{\delta(\theta-\theta_{1})-\delta(\theta+\theta_{1})}{2},
\end{eqnarray}

\noindent
where we employ the symbol $\Ev$ for the operation of taking the even part
and the symbol $\Od$ for the operation of taking the odd part, and where
we used the fact that $\delta(-\alpha)=\delta(\alpha)$, that is, the fact
that the delta ``function'' centered at zero is even. Therefore, by
representing by means of inner analytic functions within the open unit
disk the even and odd singular ``functions'' on the right-hand sides of
the last two equations above, each with two singularities, located at
$\theta_{1}$ and at $-\theta_{1}$ on the unit circle, one may recover the
representation of the delta ``function'' with a single arbitrary singular
point $\theta_{1}$.

One may restrict the general definition given above in a more significant
way by restricting in another way the set of inner analytic functions to
be considered. A restriction to a smaller set of inner analytic functions,
which we will adopt here as sufficient for our current purposes, is to
consider only the set of inner analytic functions whose sequences of
Taylor-Fourier coefficients $a_{k}$, as defined in~\cite{FTotCPI}, do not
increase with $k$ as $k\to\infty$ faster than all powers of $k$. In more
precise terms, we choose to restrict our set of inner analytic functions
to those such that, given the sequence of Taylor-Fourier coefficients
$a_{k}$, there is an integer $p>0$ such that

\begin{equation}\label{limitcond1}
  \lim_{k\to\infty}
  \frac{a_{k}}{k^{p}}
  =
  0.
\end{equation}

\noindent
Note that, since the sequence of Fourier coefficients $a_{k}$ of any
integrable real function is necessarily limited, this restriction does not
exclude any such real functions. This condition defines a definite set of
inner analytic functions within the open unit disk, and hence a definite
set of generalized functions on the unit circle. This is a sufficiently
general set of generalized functions for our purposes here, since, as we
will see, it includes the Dirac delta ``function'' and its derivatives of
arbitrarily high orders, as well as a certain class of non-integrable real
functions.

As was discussed in the previous papers~\cite{FTotCPI} and~\cite{FTotCPIV}
mentioned above, besides DP real functions and their corresponding inner
analytic functions it may also be useful to consider rotated inner
analytic functions. These are analytic functions that correspond to
shifted real functions such as $f(\theta-\theta_{1})$ and which reduce to
a real function on a diameter of the unit disk forming an angle
$\theta_{1}$ with the real axis. They constitute a fairly simple
generalization of the original structure. We will use this generalization
when we discuss a set of singular generalized functions in
Section~\ref{SECgenfunc}. However, in order to prepare the ground for the
discussion of generalized functions, we will first consider the complete
and detailed generalization of the correspondence between real functions
and complex analytic functions to all integrable real functions,
regardless of any parity considerations. The case of the rotated inner
analytic functions will then become a particular case of this larger
generalization.

Let us conclude this introduction with a note about the concept of
integrability of real functions. What we mean by integrability of real
functions in this paper is integrability in the sense of Lebesgue, with
the use of the usual Lebesgue measure. We will assume that all the real
functions under discussion here are measurable in this Lebesgue measure,
regardless of whether or not they are integrable on their whole domain.
Therefore whenever we speak of real functions in this paper, it should be
understood that we mean Lebesgue-measurable real functions. We will then
use the following result from the theory of measure and integration:
within the set of all Lebesgue-measurable real functions defined on a
compact interval, the conditions of integrability and of absolute
integrability are two equivalent conditions~\cite{RealAnalysis}. Therefore
we will use the concepts of integrability and of absolute integrability
interchangeably, as convenience requires. An alternative integrability
condition over the real functions, that we will also use, is the
requirement that they be integrable in all closed sub-intervals of
$[-\pi,\pi]$, which we will refer to as the condition of ``local
integrability''. It can be shown that in this context this condition is in
fact equivalent to the other two conditions of integrability. A simple
proof of this equivalence can be found in Appendix~\ref{APPequivinteg}.

\section{The General Case for Integrable Real Functions}\label{SECgenreal}

Let us discuss the generalization of our results for DP real functions,
obtained in~\cite{FTotCPI} and~\cite{FTotCPIV}, to a more general set of
real functions. Let it be understood that, when we talk of real functions
in this section, we always mean integrable zero-average real functions,
but not necessarily DP real functions. As was shown in~\cite{FTotCPIV},
any integrable DP real function $f(\theta)$ is representable almost
everywhere by its Fourier coefficients $a_{k}$ and is associated to an
inner analytic function $w(z)$, whose Taylor coefficients are $a_{k}$, and
from which it can be recovered almost everywhere in the $\rho\to 1$ limit,
where $z=\rho\exp(\ii\theta)$. This is so because, given an even
integrable real function $f_{\rm c}(\theta)$ we have that

\noindent
\begin{eqnarray*}
  f_{\rm c}(\theta)
  & = &
  \lim_{\rho\to 1}
  \Re[w(z)]
  \\
  & = &
  \lim_{\rho\to 1}
  \Re\!
  \left[
    f_{\rm c}(\rho,\theta)+\ii\bar{f}_{\rm c}(\rho,\theta)
  \right]
  \\
  & = &
  \lim_{\rho\to 1}
  f_{\rm c}(\rho,\theta),
\end{eqnarray*}

\noindent
where $\bar{f}_{\rm c}(\rho,\theta)$ is the Fourier-Conjugate (FC)
function to $f_{\rm c}(\rho,\theta)$, as defined in~\cite{FTotCPI}, and
where the inner analytic function associated to $f_{\rm c}(\theta)$ is
given by

\begin{displaymath}
  w(z)
  =
  f_{\rm c}(\rho,\theta)+\ii\bar{f}_{\rm c}(\rho,\theta).
\end{displaymath}

\noindent
Similarly, given an odd integrable real function $f_{\rm s}(\theta)$ we
have that

\noindent
\begin{eqnarray*}
  f_{\rm s}(\theta)
  & = &
  \lim_{\rho\to 1}
  \Im[w(z)]
  \\
  & = &
  \lim_{\rho\to 1}
  \Im\!
  \left[
    \bar{f}_{\rm s}(\rho,\theta)+\ii f_{\rm s}(\rho,\theta)
  \right]
  \\
  & = &
  \lim_{\rho\to 1}
  f_{\rm s}(\rho,\theta),
\end{eqnarray*}

\noindent
where $\bar{f}_{\rm s}(\rho,\theta)$ is the FC function to $f_{\rm
  s}(\rho,\theta)$ and where the inner analytic function associated to
$f_{\rm s}(\theta)$ is given by

\begin{displaymath}
  w(z)
  =
  \bar{f}_{\rm s}(\rho,\theta)+\ii f_{\rm s}(\rho,\theta).
\end{displaymath}

\noindent
It is immediately apparent that, if we add two DP real functions with the
same parity, we must simply add the two corresponding inner analytic
functions in order to get the inner analytic function associated to the
sum. For example, given two even DP real functions $f_{1,\rm c}(\theta)$
and $f_{2,\rm c}(\theta)$, corresponding respectively to the inner
analytic functions $w_{1}(z)$ and $w_{2}(z)$, we have that

\noindent
\begin{eqnarray*}
  f_{\rm c}(\theta)
  & = &
  f_{1,\rm c}(\theta)+f_{2,\rm c}(\theta)
  \\
  & = &
  \lim_{\rho\to 1}
  \Re\!
  \left[
    w_{1}(z)+w_{2}(z)
  \right],
\end{eqnarray*}

\noindent
so that the inner analytic function corresponding to the sum $f_{\rm
  c}(\theta)$ is simply the sum

\begin{displaymath}
  w(z)
  =
  w_{1}(z)+w_{2}(z).
\end{displaymath}

\noindent
The same is true for the addition of two odd DP real functions $f_{1,\rm
  s}(\theta)$ and $f_{2,\rm s}(\theta)$, and the corresponding inner
analytic functions $w_{1}(z)$ and $w_{2}(z)$, since in this case we have
that

\noindent
\begin{eqnarray*}
  f_{\rm s}(\theta)
  & = &
  f_{1,\rm s}(\theta)+f_{2,\rm s}(\theta)
  \\
  & = &
  \lim_{\rho\to 1}
  \Im\!
  \left[
    w_{1}(z)+w_{2}(z)
  \right],
\end{eqnarray*}

\noindent
so that once more we get for the inner analytic function associated to the
sum $f_{\rm s}(\theta)$ simply the sum

\begin{displaymath}
  w(z)
  =
  w_{1}(z)+w_{2}(z).
\end{displaymath}

\noindent
However, if we add two DP real functions with opposite parities, say
$f_{1,\rm c}(\theta)$ and $f_{2,\rm s}(\theta)$, then the situation
changes a little, since in this case we have

\noindent
\begin{eqnarray*}
  f_{1,\rm c}(\theta)
  & = &
  \lim_{\rho\to 1}
  \Re[w_{1}(z)]
  \\
  & = &
  \lim_{\rho\to 1}
  \Re\!
  \left[
    f_{1,\rm c}(\rho,\theta)+\ii\bar{f}_{1,\rm c}(\rho,\theta)
  \right],
  \\
  f_{2,\rm s}(\theta)
  & = &
  \lim_{\rho\to 1}
  \Im[w_{2}(z)]
  \\
  & = &
  \lim_{\rho\to 1}
  \Im\!
  \left[
    \bar{f}_{2,\rm s}(\rho,\theta)+\ii f_{2,\rm s}(\rho,\theta)
  \right].
\end{eqnarray*}

\noindent
In this case we can construct an inner analytic function $w(z)$ such that
the sum $f(\theta)$ of the two real functions is obtained from the
$\rho\to 1$ limit of the real part of $w(z)$ by making a complex linear
combination of the two corresponding inner analytic functions,

\noindent
\begin{eqnarray*}
  f(\theta)
  & = &
  f_{1,\rm c}(\theta)+f_{2,\rm s}(\theta)
  \\
  & = &
  \lim_{\rho\to 1}
  \Re[w_{1}(z)-\ii w_{2}(z)]
  \\
  & = &
  \lim_{\rho\to 1}
  \Re\!
  \left\{
    \left[
      f_{1,\rm c}(\rho,\theta)
      +
      f_{2,\rm s}(\rho,\theta)
    \right]
    +
    \ii
    \left[
      \bar{f}_{1,\rm c}(\rho,\theta)
      -
      \bar{f}_{2,\rm s}(\rho,\theta)
    \right]
  \right\}
  \\
  & = &
  \lim_{\rho\to 1}
  \left[
    f_{1,\rm c}(\rho,\theta)
    +
    f_{2,\rm s}(\rho,\theta)
  \right]
  \\
  & = &
  f_{\rm c}(\theta)
  +
  f_{\rm s}(\theta),
\end{eqnarray*}

\noindent
where $f_{\rm c}(\theta)=f_{1,\rm c}(\theta)$ is the even part of
$f(\theta)$ and $f_{\rm s}(\theta)=f_{2,\rm s}(\theta)$ is its odd part.
Since any real function $f(\theta)$ can be separated into its unique even
and odd parts, this gives us an inner analytic function $w(z)$ from the
real part of which this arbitrary real function can be obtained in the
$\rho\to 1$ limit. This inner analytic function is given by

\begin{equation}\label{newinner}
  w(z)
  =
  w_{1}(z)-\ii w_{2}(z).
\end{equation}

\noindent
This solution for the inner analytic function corresponding to $f(\theta)$
is not unique, in the sense that it is possible to define another one,
from the imaginary part of which the real function can also be obtained in
the $\rho\to 1$ limit,

\noindent
\begin{eqnarray*}
  f(\theta)
  & = &
  f_{1,\rm c}(\theta)+f_{2,\rm s}(\theta)
  \\
  & = &
  \lim_{\rho\to 1}
  \Im[\ii w_{1}(z)+w_{2}(z)]
  \\
  & = &
  \lim_{\rho\to 1}
  \Im\!
  \left\{
    \left[
      -
      \bar{f}_{1,\rm c}(\rho,\theta)
      +
      \bar{f}_{2,\rm s}(\rho,\theta)
    \right]
    +
    \ii
    \left[
      f_{1,\rm c}(\rho,\theta)
      +
      f_{2,\rm s}(\rho,\theta)
    \right]
  \right\}
  \\
  & = &
  \lim_{\rho\to 1}
  \left[
    f_{1,\rm c}(\rho,\theta)
    +
    f_{2,\rm s}(\rho,\theta)
  \right]
  \\
  & = &
  f_{\rm c}(\theta)
  +
  f_{\rm s}(\theta),
\end{eqnarray*}

\noindent
so that we have for this alternative inner analytic function

\begin{displaymath}
  w'(z)
  =
  \ii w_{1}(z)+w_{2}(z).
\end{displaymath}

\noindent
However, this second solution need not concern us, since it is in fact
proportional to the first one, for we can see that $w'(z)=\ii w(z)$. So
long as we agree that the real functions are to be obtained as the
$\rho\to 1$ limits of, say, the real parts of the inner analytic
functions, there is a unique inner analytic function associated to each
real function. This can be applied to the purely odd real functions as
well, if we associate to them the inner analytic function $-\ii w(z)$,
rather than $w(z)$ as we have been doing so far, since we have that

\begin{displaymath}
  \Im[w(z)]
  =
  \Re[-\ii w(z)].
\end{displaymath}

\noindent
In short, given any real function $f(\theta)$, regardless of whether or
not it has definite parity properties, we know how to build from it the
unique inner analytic function that gives us back that real function as
the $\rho\to 1$ limit of its real part: first we separate the real
function $f(\theta)$ into its even and odd parts $f_{1,\rm c}(\theta)$ and
$f_{2,\rm s}(\theta)$; then, we use the previous method of definition, as
given in~\cite{FTotCPI}, to determine the corresponding old-style inner
analytic functions $w_{1}(z)$ and $w_{2}(z)$; finally, we define the inner
analytic function corresponding to $f(\theta)$ as

\begin{displaymath}
  w(z)
  =
  w_{1}(z)
  -
  \ii
  w_{2}(z).
\end{displaymath}

\noindent
Let us now consider the inverse problem, that is, given an arbitrary real
function $f(\theta)$ and its corresponding inner analytic function $w(z)$,
defined according to our new criterion, the problem of how to obtain from
it the inner analytic functions corresponding to the even and odd parts of
the real function. We may assume that the inner analytic function is given
as

\begin{displaymath}
  w(z)
  =
  f(\rho,\theta)
  +
  \ii
  \bar{f}(\rho,\theta),
\end{displaymath}

\noindent
where $f(\rho,\theta)$ is harmonic and $\bar{f}(\rho,\theta)$ is the
harmonic conjugate function to $f(\rho,\theta)$, and where we assume
adherence to the criterion that the real function is to be recovered from
the $\rho\to 1$ limit of the real part of the corresponding inner analytic
function,

\noindent
\begin{eqnarray*}
  f(\theta)
  & = &
  \lim_{\rho\to 1}
  \Re[w(z)]
  \\
  & = &
  \lim_{\rho\to 1}
  \Re[f(\rho,\theta)+\ii\bar{f}(\rho,\theta)]
  \\
  & = &
  \lim_{\rho\to 1}
  f(\rho,\theta).
\end{eqnarray*}

\noindent
On the other hand, we assume that $f(\theta)$ is decomposed into its even
and odd parts as

\begin{displaymath}
  f(\theta)
  =
  f_{1,\rm c}(\theta)+f_{2,\rm s}(\theta),
\end{displaymath}

\noindent
where the two DP real functions $f_{1,\rm c}(\theta)$ and $f_{2,\rm
  s}(\theta)$ are associated to the two inner analytic functions
$w_{1}(z)$ and $w_{2}(z)$,

\noindent
\begin{eqnarray*}
  f_{1,\rm c}(\theta)
  & = &
  \lim_{\rho\to 1}
  \Re[w_{1}(z)],
  \\
  f_{2,\rm s}(\theta)
  & = &
  \lim_{\rho\to 1}
  \Re[w_{2}(z)],
\end{eqnarray*}

\noindent
which are defined according to our new criterion,

\noindent
\begin{eqnarray*}
  w_{1}(z)
  & = &
  f_{1,\rm c}(\rho,\theta)+\ii\bar{f}_{1,\rm c}(\rho,\theta),
  \\
  w_{2}(z)
  & = &
  f_{2,\rm s}(\rho,\theta)-\ii\bar{f}_{2,\rm s}(\rho,\theta).
\end{eqnarray*}

\noindent
Since we have that $w(z)=w_{1}(z)+w_{2}(z)$ we now see that we have for
the inner analytic function $w(z)$ corresponding to $f(\theta)$, written
in terms of the real and imaginary parts of the inner analytic functions
$w_{1}(z)$ and $w_{2}(z)$ corresponding respectively to the even and odd
parts of $f(\theta)$,

\noindent
\begin{eqnarray*}
  \left[
    f(\rho,\theta)
    +
    \ii
    \bar{f}(\rho,\theta)
  \right]
  & = &
  \left[
    f_{1,\rm c}(\rho,\theta)
    +
    \ii
    \bar{f}_{1,\rm c}(\rho,\theta)
  \right]
  +
  \left[
    f_{2,\rm s}(\rho,\theta)
    -
    \ii
    \bar{f}_{2,\rm s}(\rho,\theta)
  \right]
  \\
  & = &
  \left[
    f_{1,\rm c}(\rho,\theta)
    +
    f_{2,\rm s}(\rho,\theta)
  \right]
  +
  \ii
  \left[
    \bar{f}_{1,\rm c}(\rho,\theta)
    -
    \bar{f}_{2,\rm s}(\rho,\theta)
  \right],
\end{eqnarray*}

\noindent
which implies that we have

\noindent
\begin{eqnarray*}
  f(\rho,\theta)
  & = &
  f_{1,\rm c}(\rho,\theta)
  +
  f_{2,\rm s}(\rho,\theta),
  \\
  \bar{f}(\rho,\theta)
  & = &
  \bar{f}_{1,\rm c}(\rho,\theta)
  -
  \bar{f}_{2,\rm s}(\rho,\theta).
\end{eqnarray*}

\noindent
These two relations among real functions in turn imply that the DP
functions on the right-hand sides are the even and odd parts with respect
to $\theta$ of the functions in the left-hand sides. Characterizing the
functions by their parities, and remembering that Fourier Conjugation
reverses parity, we have

\noindent
\begin{eqnarray*}
  f_{1,\rm c}(\rho,\theta)
  & = &
  \Ev[f(\rho,\theta)],
  \\
  f_{2,\rm s}(\rho,\theta)
  & = &
  \Od[f(\rho,\theta)],
  \\
  \bar{f}_{1,\rm c}(\rho,\theta)
  & = &
  \Od[\bar{f}(\rho,\theta)],
  \\
  \bar{f}_{2,\rm s}(\rho,\theta)
  & = &
  -
  \Ev[\bar{f}(\rho,\theta)],
\end{eqnarray*}

\noindent
where the symbols $\Ev$ and $\Od$ relate to the parity with respect to
$\theta$. As a consequence, we have for the inner analytic functions
associated to the even and odd parts of $f(\theta)$, respectively

\noindent
\begin{eqnarray*}
  w_{1}(z)
  & = &
  f_{1,\rm c}(\rho,\theta)
  +
  \ii
  \bar{f}_{1,\rm c}(\rho,\theta)
  \\
  & = &
  \Ev[f(\rho,\theta)]
  +
  \ii
  \Od[\bar{f}(\rho,\theta)],
  \\
  w_{2}(z)
  & = &
  f_{2,\rm s}(\rho,\theta)
  -
  \ii
  \bar{f}_{2,\rm s}(\rho,\theta)
  \\
  & = &
  \Od[f(\rho,\theta)]
  +
  \ii
  \Ev[\bar{f}(\rho,\theta)],
\end{eqnarray*}

\noindent
where we are still using the new standard criterion that the real
functions are to be recovered from the $\rho\to 1$ limits of the real
parts of the corresponding inner analytic functions. In conclusion, given
the inner analytic function $w(z)$ corresponding to $f(\theta)$, we may
write for the inner analytic function $w_{1}(z)$ corresponding to the even
part of $f(\theta)$, and for the inner analytic function $w_{2}(z)$
corresponding to the odd part of $f(\theta)$,

\noindent
\begin{eqnarray*}
  w_{1}(z)
  & = &
  \Ev\{\Re[w(z)]\}
  +
  \ii\,
  \Od\{\Im[w(z)]\},
  \\
  w_{2}(z)
  & = &
  \Od\{\Re[w(z)]\}
  +
  \ii\,
  \Ev\{\Im[w(z)]\},
\end{eqnarray*}

\noindent
which are both, therefore, uniquely and completely determined. We have,
therefore, a consistent way to associate arbitrary real functions with
corresponding inner analytic functions, in such a way that each real
function is recovered as the $\rho\to 1$ limit of the real part of the
corresponding inner analytic function. Due to this, the sum of any pair of
real functions is now related to the simple sum of the corresponding pair
of inner analytic functions. In other words, if we adhere to this new
standard way to relate the real functions and the inner analytic
functions, then the set of inner analytic functions is seen to inherit
from the set of real functions its character as a vector space with real
scalars.

One interesting special case, which will be of much use to us here, is
that of a rotated inner analytic function. In order to discuss it, let us
take the case of an even real function $f(\theta)$, and the corresponding
inner analytic function $w(z)$. Let us suppose that we generate from
$f(\theta)$ another function by just shifting the variable $\theta$ by a
real constant $\theta_{1}$, in order to define the shifted function
$f_{1}(\theta)=f(\theta-\theta_{1})$. Clearly this modified function is no
longer an even function of $\theta$. It has now non-zero even and odd
parts, and the results previously obtained in this section can all be
applied to it. Let us determine how the inner analytic function $w(z)$
corresponding to $f(\theta)$ changes into a new analytic function
$w_{1}(z)$ that corresponds to the shifted real function
$f_{1}(\theta)$. Let us recall that we have

\begin{displaymath}
  z
  =
  \rho\e{\ii\theta},
\end{displaymath}

\noindent
which is the only place there $\theta$ appears in the construction of the
inner analytic function $w(z)$, and which, upon the shift
$\theta\to\theta-\theta_{1}$, changes to

\begin{displaymath}
  z'
  =
  \rho\e{\ii(\theta-\theta_{1})}.
\end{displaymath}

\noindent
If we define the complex constant $z_{1}$ associated to $\theta_{1}$ as a
point on the unit circle,

\begin{displaymath}
  z_{1}
  =
  \e{\ii\theta_{1}},
\end{displaymath}

\noindent
then we have

\noindent
\begin{eqnarray*}
  z'
  & = &
  \rho\e{\ii\theta}\e{-\ii\theta_{1}}
  \\
  & = &
  \frac{z}{z_{1}}.
\end{eqnarray*}

\noindent
It follows therefore that this shift in the real variable $\theta$
corresponds to the exchange of $z$ for the ratio $z/z_{1}$ in the inner
analytic function, so that we simply have

\begin{displaymath}
  w_{1}(z)
  =
  w(z/z_{1}).
\end{displaymath}

\noindent
Since $w(z)$ reduces to a real function on the interval $(-1,1)$ of the
real axis, it follows that $w_{1}(z)$ reduces to a real function when
$z/z_{1}$ is real and within that interval. This is realized when $z$ and
$z_{1}$ are collinear in the complex plane, which means that we have
$\theta=\theta_{1}$ or $\theta=\pi+\theta_{1}$. Therefore, the inner
analytic function $w_{1}(z)$ reduces to a real function on a diameter of
the unit circle that makes an angle $\theta_{1}$ with the real axis.
Similar conclusions are reached when $f(\theta)$ is an odd function. In
this case the modified inner analytic function $w_{1}(z)$ is still given
by the formula above, but it now reduces to a purely imaginary function
over that same diagonal of the unit circle.

\subsection{Updating Concepts and Definitions}

We see therefore that the correspondence between real functions and inner
analytic functions can be generalized in a simple and straightforward way
to the set of all possible integrable real functions, regardless of any
parity considerations. Upon this generalization the new class of inner
analytic functions $w(z)$ that arises is somewhat larger than the original
one. These are all still analytic functions within the open unit disk, and
they do all have the property that $w(0)=0$, just like the original ones,
since we did not lift here the limitation that the real functions be
zero-average. However, they no longer reduce to real functions over the
interval $(-1,1)$ of the real axis, nor they do so, in general, over any
diameter of the unit disk. Only if the given real function is even will
the corresponding inner analytic function reduce to a real function on
that interval of the real axis. On the other hand, if the given real
function is odd, then the corresponding inner analytic function, if
defined according to our new criterion, reduces to a purely imaginary
function on that same interval. It should be noted that, although this new
class of inner analytic functions is a generalization of the previous one,
it is still far from including all possible analytic functions within the
open unit disk.

We see therefore that we have been led to generalize the concept of an
inner analytic function when we worked out the generalization of our
previous structure to non-DP real functions. It follows that several other
concepts and definitions are in need of revision as well. For example, it
is no longer true that the coefficients of the complex Taylor series are
real. Also, the set of Fourier coefficients of the real functions consist
now of two sequences, one for the even part and another for the odd part
of the real function. Therefore, we can no longer use the notation $a_{k}$
for the Fourier coefficients as we did before, where $a_{k}$ could be the
name for either the Fourier coefficients $\alpha_{k}$ of the cosine series
of the even part or the Fourier coefficients $\beta_{k}$ of the sine
series of the odd part, as the case may be. We have now two sequences of
Taylor-Fourier coefficients associated to each real function, and we must
revert to the traditional notation $\alpha_{k}$ and $\beta_{k}$. In fact,
given the new form of the inner analytic functions, shown in
Equation~(\ref{newinner}), it follows that the coefficients of the complex
power series are now given by

\begin{displaymath}
  c_{k}
  =
  \alpha_{k}
  -
  \ii
  \beta_{k}.
\end{displaymath}

\noindent
Since $\alpha_{k}$ and $\beta_{k}$ are both essentially arbitrary
sequences of independent real coefficients, it follows that $c_{k}$ is a
sequence of essentially arbitrary complex coefficients. However, if we
restrict ourselves to normal integrable real functions on the unit circle,
then these coefficients are not really completely arbitrary, since
$\alpha_{k}$ and $\beta_{k}$ must satisfy the condition that they can be
obtained as the usual integrals involving the real functions. Due to this,
they are both limited sequences of coefficients, and therefore so is
$c_{k}$.

Another concept that is in need of revision is that of Fourier
Conjugation. In the previous paper~\cite{FTotCPI} we defined the concept
of Fourier Conjugation by the simple interchange of cosines and sines in
the Fourier series of the DP real functions. This definition must now be
adapted to the new larger structure, where it will still be related, as it
was before, to the concept of harmonic conjugation. Let us state here our
main new definitions.

\begin{description}

\item[Inner analytic function:] an inner analytic function $w(z)$ is a
  complex function that is analytic within the open unit disk and that has
  the property that $w(0)=0$.

  Additionally, for the purposes of this paper we will adopt the further
  restriction that the Taylor-Fourier coefficients $c_{k}$ of $w(z)$
  around the origin satisfy the property that there is an integer $p>0$
  such that

  \begin{displaymath}
    \lim_{k\to\infty}
    \frac{|c_{k}|}{k^{p}}
    =
    0.
  \end{displaymath}

\item[Taylor-Fourier Coefficients:] given a zero-average integrable real
  function $f(\theta)$, one obtains from it the two sequences of Fourier
  coefficients, $\alpha_{k}$ from its even part and $\beta_{k}$ from its
  odd part, in either case for $k\geq 1$. The Taylor-Fourier coefficients
  are then defined as

  \begin{displaymath}
    c_{k}
    =
    \alpha_{k}
    -
    \ii
    \beta_{k},
  \end{displaymath}

  \noindent
  so that we have the inner analytic function within the open unit disk,

  \begin{displaymath}
    w(z)
    =
    \sum_{k=1}^{\infty}
    c_{k}z^{k},
  \end{displaymath}

  \noindent
  from the real part of which $f(\theta)$ is recovered in the $\rho\to 1$
  limit to the unit circle.

\item[Fourier Conjugate:] given that the inner analytic function $w(z)$
  may be written within the open unit disk as

  \begin{displaymath}
    w(z)
    =
    f(\rho,\theta)
    +
    \ii
    \bar{f}(\rho,\theta),
  \end{displaymath}

  \noindent
  where $f(\rho,\theta)$ is harmonic and $\bar{f}(\rho,\theta)$ is the
  harmonic conjugate function to $f(\rho,\theta)$, and that the
  corresponding real function $f(\theta)$ is given by

  \begin{displaymath}
    f(\theta)
    =
    \lim_{\rho\to 1}
    f(\rho,\theta),
  \end{displaymath}

  \noindent
  we define the Fourier Conjugate real function $\bar{f}(\theta)$ to the
  real function $f(\theta)$ as

  \begin{displaymath}
    \bar{f}(\theta)
    =
    \lim_{\rho\to 1}
    \bar{f}(\rho,\theta).
  \end{displaymath}

\end{description}

\noindent
The largest possible generalization of our structure is obtained if we
enlarge the set of inner analytic functions to all complex analytic
function within the open unit disk that have the property that $w(0)=0$.
The corresponding set of complex power series includes all such series
that converge everywhere on the open unit disk and that have no constant
term. There is then a corresponding generalization of the real objects on
the unit circle. The new set includes not only all integrable real
functions, but also singular objects we will call generalized functions,
and some non-integrable real functions as well. As stated in the
introduction, as a practical measure of simplification we will limit our
set of inner analytic functions, for the purposes of this paper, to those
whose Taylor-Fourier coefficients $c_{k}$ satisfy the condition that there
is an integer $p>0$ such that

\begin{equation}\label{limitcond2}
  \lim_{k\to\infty}
  \frac{|c_{k}|}{k^{p}}
  =
  0.
\end{equation}

\noindent
In the next two sections we will push this generalization forward,
enlarging the whole structure to include a set of singular generalized
functions, as well as a class of singular non-integrable functions. As we
will see, all the operations discussed here in the context of the
integrable real functions can be applied to these generalized cases as
well. This is due to the fact that both the normal functions and the
generalized functions can be represented by inner analytic functions, so
that all operations performed on either normal or generalized functions
can be interpreted in terms of corresponding operations on the inner
analytic functions.

\section{A Set of Singular Generalized Functions}\label{SECgenfunc}

Let us start by reviewing the analysis of the representation within the
open unit disk of the Dirac delta ``function'' that was presented in
detail in~\cite{FTotCPI}. Consider the very simple analytic function of
$z$

\begin{equation}\label{innerdelta}
  w_{\delta}(z,z_{1})
  =
  \frac{1}{2\pi}
  -
  \frac{1}{\pi}\,
  \frac{z}{z-z_{1}},
\end{equation}

\noindent
where $z=\rho\exp(\ii\theta)$ and $z_{1}=\exp(\ii\theta_{1})$, which
corresponds to the sequence of Taylor-Fourier coefficients
$\alpha_{k}=1/\pi$ for all $k\geq 1$, when expanded in powers of
$z/z_{1}$, as one can see by the straightforward calculation of its Taylor
coefficients, or by simply expanding the ratio shown in the form of a
geometric series,

\begin{equation}\label{expandelta}
  w_{\delta}(z,z_{1})
  =
  \frac{1}{2\pi}
  +
  \frac{1}{\pi}
  \sum_{k=1}^{\infty}
  \left(\frac{z}{z_{1}}\right)^{k}.
\end{equation}

\noindent
Note that the coefficients $\alpha_{k}$ correspond to a Fourier series
which diverges almost everywhere. The analytic function of $z$ shown in
Equations~(\ref{innerdelta}) and~(\ref{expandelta}) is an extended inner
analytic function, that is, a rotated inner analytic function with the
constant shown added, and its real part represents the Dirac delta
``function'' $\delta(\theta-\theta_{1})$ when one takes the limit to the
unit circle,

\noindent
\begin{eqnarray*}
  \delta(\theta-\theta_{1})
  & = &
  \lim_{\rho\to 1}
  \Re[w_{\delta}(z,z_{1})]
  \\
  & = &
  \frac{1}{2\pi}
  -
  \frac{1}{\pi}\,
  \lim_{\rho\to 1}
  \frac
  {\rho^{2}-\rho\cos(\Delta\theta)}
  {\left(\rho^{2}+1\right)-2\rho\cos(\Delta\theta)},
\end{eqnarray*}

\noindent
where one can see that the dependence is indeed only on the difference
$\Delta\theta=\theta-\theta_{1}$ and that the ``function'' is even on
$\Delta\theta$. The limit represents the delta ``function'' in the sense
that it has the following properties, as demonstrated in~\cite{FTotCPI}:

\begin{enumerate}

\item The definition of $\delta(\theta-\theta_{1})$ tends to zero when one
  takes the $\rho\to 1$ limit while keeping $\theta\neq\theta_{1}$;

\item The definition of $\delta(\theta-\theta_{1})$ diverges to positive
  infinity when one takes the $\rho\to 1$ limit with $\theta=\theta_{1}$;

\item In the $\rho\to 1$ limit the integral

  \begin{displaymath}
    \int_{a}^{b}d\theta\,
    \delta(\theta-\theta_{1})
    =
    1,
  \end{displaymath}

  \noindent
  has the value shown, for any interval $[a,b]$ which contains the point
  $\theta_{1}$;

\item Given any continuous function $g(\theta)$, in the $\rho\to 1$ limit
  the integral

  \begin{displaymath}
    \int_{a}^{b}d\theta\,
    g(\theta)
    \delta(\theta-\theta_{1})
    =
    g(\theta_{1}),
  \end{displaymath}

  \noindent
  has the value shown, for any interval $[a,b]$ which contains the point
  $\theta_{1}$.

\end{enumerate}

\noindent
In fact, as shown in~\cite{FTotCPI}, in this last property the hypothesis
that $g(\theta)$ is continuous on $\theta$ can be relaxed, since
$g(\theta)$, being a real function of $\theta$, is obtained as the
$\rho\to 1$ limit of a harmonic function $g(\rho,\theta)$, and it is
sufficient that this last function be continuous in $\rho$ as we make
$\rho\to 1$. Also note that, although it is customary to list both
separately, the third property is in fact just a particular case of the
fourth property. The latter is the most important one, since it
establishes the delta ``function'' as the integration kernel of a
real-valued linear functional acting in the space of integrable real
functions, defined by the integral of the product shown. This is the type
of interpretation of this singular object that is established in the
Schwartz theory of distributions~\cite{DistTheory}.

\subsection{Even and Odd Parts of the Delta ``Function''}

Let us take the opportunity to use the extended inner analytic function in
Equation~(\ref{innerdelta}), as well as the structure developed in
Section~\ref{SECgenreal}, to exemplify the determination of the extended
inner analytic functions associated to the even and odd parts of the delta
``function''. This is, in essence, the definition of two new extended
inner analytic functions and of the corresponding real singular
generalized functions. In order to do this we must first determine the
extended inner analytic function associated to the reflected delta
``function'', which is given by

\begin{displaymath}
  \delta(-\theta-\theta_{1})
  =
  \delta(\theta+\theta_{1}),
\end{displaymath}

\noindent
since the delta ``function'' centered at zero is even. We are now in a
position to use the definitions in Equation~(\ref{deltaparity}) for the
even and odd parts of the delta ``function'' $\delta(\theta-\theta_{1})$.
Starting from the definition of the delta function with the usual
argument,

\noindent
\begin{eqnarray*}
  \delta(\theta-\theta_{1})
  & = &
  \lim_{\rho\to 1}
  \Re[w_{\delta}(z,z_{1})],
  \\
  w_{\delta}(z,z_{1})
  & = &
  \frac{1}{2\pi}
  -
  \frac{1}{\pi}\,
  \frac{z}{z-z_{1}},
\end{eqnarray*}

\noindent
we may now write the definition of the same function with the sign of
$\theta_{1}$ reversed. It is easy to see that simply exchanging the sign
of $\theta_{1}$ has the effect that $z_{1}=\exp(\ii\theta_{1})$ gets
mapped onto its complex conjugate, so that we have

\noindent
\begin{eqnarray*}
  \delta(\theta+\theta_{1})
  & = &
  \lim_{\rho\to 1}
  \Re\!\left[w_{\delta}\left(z,z_{1}^{*}\right)\right],
  \\
  w_{\delta}\!\left(z,z_{1}^{*}\right)
  & = &
  \frac{1}{2\pi}
  -
  \frac{1}{\pi}\,
  \frac{z}{z-z_{1}^{*}}.
\end{eqnarray*}

\noindent
Observe that the complex conjugation applies only to the constant $z_{1}$,
so that $w_{\delta}\!\left(z,z_{1}^{*}\right)$ is still an analytic
function of $z$. From the definitions in Equation~(\ref{deltaparity}) and
the fact that we may simply add the corresponding inner analytic functions
within the open unit disk, it follows that the extended inner analytic
function $w_{1}(z)$ corresponding to the even part and the inner analytic
function $w_{2}(z)$ corresponding to the odd part are given by

\noindent
\begin{eqnarray*}
  w_{1}(z)
  & = &
  \frac{1}{2\pi}
  -
  \frac{z}{2\pi}
  \left(
    \frac{1}{z-z_{1}}
    +
    \frac{1}{z-z_{1}^{*}}
  \right)
  \\
  & = &
  \frac{1}{2\pi}
  -
  \frac{z}{2\pi}\,
  \frac
  {2z-\left(z_{1}+z_{1}^{*}\right)}
  {(z-z_{1})\left(z-z_{1}^{*}\right)},
  \\
  w_{2}(z)
  & = &
  -\,
  \frac{z}{2\pi}
  \left(
    \frac{1}{z-z_{1}}
    -
    \frac{1}{z-z_{1}^{*}}
  \right)
  \\
  & = &
  -\,
  \frac{z}{2\pi}\,
  \frac
  {\left(z_{1}-z_{1}^{*}\right)}
  {(z-z_{1})\left(z-z_{1}^{*}\right)}.
\end{eqnarray*}

\noindent
Note that these are not rotated inner analytic functions, but old-style
ones except for the constant term of $w_{1}(z)$ and the lack of a factor
of $\ii$ in $w_{2}(z)$, so that we have that $w_{1}(0)=1/(2\pi)$ and
$w_{2}(0)=0$, and also that $w_{1}(z)$ reduces to a real function on the
interval $(-1,1)$ of the real axis, while $w_{2}(z)$ reduces to a purely
imaginary function on that same interval. Therefore, taking the limit to
the unit circle, we may conclude that the even and odd parts of the delta
``function'' are represented by

\noindent
\begin{eqnarray*}
  \Ev[\delta(\theta-\theta_{1})]
  & = &
  \frac{\delta(\theta-\theta_{1})+\delta(\theta+\theta_{1})}{2}
  \\
  & = &
  \frac{1}{2\pi}
  -
  \frac{1}{2\pi}
  \lim_{\rho\to 1}
  \Re\!
  \left[
    z\,
    \frac
    {2z-\left(z_{1}+z_{1}^{*}\right)}
    {(z-z_{1})\left(z-z_{1}^{*}\right)}
  \right],
  \\
  \Od[\delta(\theta-\theta_{1})]
  & = &
  \frac{\delta(\theta-\theta_{1})-\delta(\theta+\theta_{1})}{2}
  \\
  & = &
  -\,
  \frac{1}{2\pi}
  \lim_{\rho\to 1}
  \Re\!
  \left[
    z\,
    \frac
    {\left(z_{1}-z_{1}^{*}\right)}
    {(z-z_{1})\left(z-z_{1}^{*}\right)}
  \right].
\end{eqnarray*}

\noindent
Due to the analyticity of the functions involved, these limits exist
everywhere except at the two points of singularity $z_{1}$ and
$z_{1}^{*}$. This illustrates how one can use the inner analytic functions
to define generalized functions on the unit circle, in a very simple way,
making use of the vector-space properties of the inner analytic functions,
when defined according to our new standard convention.

\subsection{Derivatives of the Delta ``Function''}

Starting once gain from the extended inner analytic function associated to
the delta ``function'' $\delta(\theta-\theta_{1})$, shown in
Equation~(\ref{innerdelta}), we may now define from it further generalized
functions with the use of the operation of integration by parts. In order
to do this we write $w_{\delta}(z,z_{1})$ as

\begin{displaymath}
  w_{\delta}(z,z_{1})
  =
  \frac{1}{2\pi}
  +
  f(\rho,\theta,\theta_{1})
  +
  \ii
  \bar{f}(\rho,\theta,\theta_{1}),
\end{displaymath}

\noindent
and consider integrals on circles of radius $\rho<1$ centered at the
origin involving a function $g(\rho,\theta)$. Since this will always be
the real part of the inner analytic function corresponding to some real
function $g(\theta)$, it follows that within the open unit disk it is
infinitely differentiable on both arguments. Observe that the following
integrals are {\em not} integrals of an analytic function on a closed
contour, but rather a pair of real integrals of real functions on the
circles, one integral in the real part and one in the imaginary part,

\noindent
\begin{eqnarray}\label{firstform1}
  I_{\delta}(\rho,\theta_{1})
  & = &
  \int_{-\pi}^{\pi}d\theta\,
  \frac{dg(\rho,\theta)}{d\theta}\,
  w_{\delta}(z,z_{1})
  \nonumber\\
  & = &
  \int_{-\pi}^{\pi}d\theta\,
  \frac{dg(\rho,\theta)}{d\theta}
  \left[
    \frac{1}{2\pi}
    +
    f(\rho,\theta,\theta_{1})
    +
    \ii
    \bar{f}(\rho,\theta,\theta_{1})
  \right].
\end{eqnarray}

\noindent
Since on the circles we have $dz=\ii zd\theta$, these integrals within the
open unit disk can also be written as

\begin{displaymath}
  I_{\delta}(\rho,\theta_{1})
  =
  \oint_{C}dz\,
  \frac{dg(\rho,\theta)}{dz}\,
  w_{\delta}(z,z_{1}),
\end{displaymath}

\noindent
where $C$ is any circle with radius $\rho<1$ centered at the origin and
where $g(\rho,\theta)$ is a {\em real} function of $\rho$ and $\theta$.
In either case, since the circle has no boundary, it is clear that we may
integrate by parts without generating an integrated term, and thus obtain

\noindent
\begin{eqnarray}\label{secondform1}
  I_{\delta}(\rho,\theta_{1})
  & = &
  -
  \oint_{C}dz\,
  g(\rho,\theta)\,
  \frac{d}{dz}
  w_{\delta}(z,z_{1})
  \nonumber\\
  & = &
  -
  \int_{-\pi}^{\pi}d\theta\,
  g(\rho,\theta)
  \left[
    \frac{d}{d\theta}
    f(\rho,\theta,\theta_{1})
    +
    \ii\,
    \frac{d}{d\theta}
    \bar{f}(\rho,\theta,\theta_{1})
  \right].
\end{eqnarray}

\noindent
Note that the integration by parts can be written in either one of these
two equivalent ways, within the open unit disk,

\noindent
\begin{eqnarray*}
  \int_{-\pi}^{\pi}d\theta\,
  \frac{dg(\rho,\theta)}{d\theta}\,
  w_{\delta}(\rho,\theta,\theta_{1})
  & = &
  -
  \int_{-\pi}^{\pi}d\theta\,
  g(\rho,\theta)\,
  \frac{dw_{\delta}(\rho,\theta,\theta_{1})}{d\theta},
  \\
  \oint_{C}dz\,
  \frac{dg(\rho,\theta)}{dz}\,
  w_{\delta}(z,z_{1})
  & = &
  -
  \oint_{C}dz\,
  g(\rho,\theta)\,
  \frac{dw_{\delta}(z,z_{1})}{dz},
\end{eqnarray*}

\noindent
where $C$ is any circle with radius $\rho<1$ centered at the origin and
where $w_{\delta}(z,z_{1})$ can be written either in terms of $z$ or in
terms of $\rho$ and $\theta$. We now observe that the real part of the
original form of this integral in Equation~(\ref{firstform1}), when the
limit $\rho\to 1$ to the unit circle is taken, becomes the integral of the
product of the delta ``function'' and the derivative of the real function,
so that we have

\noindent
\begin{eqnarray*}
  I_{\delta}(1,\theta_{1})
  & = &
  \int_{-\pi}^{\pi}d\theta\,
  \left[
    \frac{dg(1,\theta)}{d\theta}\,
    \delta(\theta-\theta_{1})
    +
    \ii\,
    \frac{dg(1,\theta)}{d\theta}\,
    \bar{f}(1,\theta,\theta_{1})
  \right]
  \\
  & = &
  \frac{dg}{d\theta}(1,\theta_{1})
  +
  \ii
  \int_{-\pi}^{\pi}d\theta\,
  \frac{dg(1,\theta)}{d\theta}\,
  \bar{f}(1,\theta,\theta_{1}),
\end{eqnarray*}

\noindent
where we used the integration properties of $\delta(\theta-\theta_{1})$.
If we also take the limit to the unit circle of the last expression in
Equation~(\ref{secondform1}), then it follows that on that circle we have

\noindent
\begin{eqnarray*}
  I_{\delta}(1,\theta_{1})
  & = &
  -
  \lim_{\rho\to 1}
  \int_{-\pi}^{\pi}d\theta\,
  g(\rho,\theta)
  \left[
    \frac{d}{d\theta}
    f(\rho,\theta,\theta_{1})
    +
    \ii\,
    \frac{d}{d\theta}
    \bar{f}(\rho,\theta,\theta_{1})
  \right]
  \\
  & = &
  -
  \int_{-\pi}^{\pi}d\theta\,
  g(1,\theta)
  \lim_{\rho\to 1}
  \left[
    \frac{d}{d\theta}
    f(\rho,\theta,\theta_{1})
  \right]
  -
  \ii
  \int_{-\pi}^{\pi}d\theta\,
  g(1,\theta)\,
  \lim_{\rho\to 1}
  \left[
    \frac{d}{d\theta}
    \bar{f}(\rho,\theta,\theta_{1})
  \right].
\end{eqnarray*}

\noindent
From these two expressions for $I_{\delta}(1,\theta_{1})$ we conclude that
we have

\noindent
\begin{eqnarray*}
  \lefteqn
  {
    \frac{dg}{d\theta}(1,\theta_{1})
    +
    \ii
    \int_{-\pi}^{\pi}d\theta\,
    \frac{dg(1,\theta)}{d\theta}\,
    \bar{f}(1,\theta,\theta_{1})
  }
  &   &
  \\
  & = &
  -
  \int_{-\pi}^{\pi}d\theta\,
  g(1,\theta)
  \lim_{\rho\to 1}
  \left[
    \frac{d}{d\theta}
    f(\rho,\theta,\theta_{1})
  \right]
  -
  \ii
  \int_{-\pi}^{\pi}d\theta\,
  g(1,\theta)
  \lim_{\rho\to 1}
  \left[
    \frac{d}{d\theta}
    \bar{f}(\rho,\theta,\theta_{1})
  \right].
\end{eqnarray*}

\noindent
The real part of this equality tell us that the expression involving the
limit of the derivative shown within the integral at the right-hand side
can be interpreted as a definition of the first derivative of the delta
``function'',

\begin{displaymath}
  \frac{dg}{d\theta}(1,\theta_{1})
  =
  -
  \int_{-\pi}^{\pi}d\theta\,
  g(1,\theta)
  \lim_{\rho\to 1}
  \left[
    \frac{d}{d\theta}
    f(\rho,\theta,\theta_{1})
  \right]
  \;\;\;\Rightarrow
\end{displaymath}
\begin{displaymath}
  \frac{d}{d\theta}
  \delta(\theta-\theta_{1})
  \equiv
  \lim_{\rho\to 1}
  \left[
    \frac{d}{d\theta}
    f(\rho,\theta,\theta_{1})
  \right].
\end{displaymath}

\noindent
This is a new generalized function, that has the property that through
integration over the unit circle it attributes to a given differentiable
function minus the value of its derivative at the point $\theta_{1}$,

\begin{displaymath}
  \int_{-\pi}^{\pi}d\theta\,
  g(1,\theta)\,
  \frac{d}{d\theta}
  \delta(\theta-\theta_{1})
  =
  -\,
  \frac{dg}{d\theta}(1,\theta_{1}).
\end{displaymath}

\noindent
It is, therefore, the integration kernel of a real-valued linear
functional acting in the space of integrable real functions. As we will
see shortly, this new generalized function is indeed dependent only on the
difference $\theta-\theta_{1}$. We will adopt for this new singular object
the notation

\begin{displaymath}
  \delta'(\theta-\theta_{1})
  \equiv
  \frac{d}{d\theta}
  \delta(\theta-\theta_{1}).
\end{displaymath}

\noindent
We may now determine the inner analytic function associated to this new
singular object. Since derivatives with respect to $\theta$ may be written
as logarithmic derivatives of the corresponding inner analytic function
within the open unit disk, as shown in~\cite{FTotCPII}, we have, using the
notation established in that paper,

\noindent
\begin{eqnarray*}
  w_{\delta'}(z,z_{1})
  & = &
  w^{\ldot}_{\delta}(z,z_{1})
  \\
  & = &
  \ii z\,
  \frac{d}{dz}
  w_{\delta}(z,z_{1})
  \\
  & = &
  \frac{d}{d\theta}
  f(\rho,\theta,\theta_{1})
  +
  \ii
  \frac{d}{d\theta}
  \bar{f}(\rho,\theta,\theta_{1}),
\end{eqnarray*}

\noindent
so that the generalized function $\delta'(\theta-\theta_{1})$ is given by
the limit to the unit circle of the real part of the inner analytic
function $w_{\delta'}(z,z_{1})$,

\begin{displaymath}
  \delta'(\theta-\theta_{1})
  =
  \lim_{\rho\to 1}
  \Re[w_{\delta'}(z,z_{1})].
\end{displaymath}

\noindent
Since we have the inner analytic function that corresponds to the delta
``function'' in explicit form, we may perform an explicit calculation in
order to obtain the inner analytic function that corresponds to the
derivative of the delta ``function'',

\noindent
\begin{eqnarray*}
  \ii z\,
  \frac{d}{dz}
  w_{\delta}(z,z_{1})
  & = &
  \ii z\,
  \frac{d}{dz}
  \left[
    \frac{1}{2\pi}
    -
    \frac{1}{\pi}\,
    \frac{z}{z-z_{1}}
  \right]
  \\
  & = &
  \frac{-\ii z}{\pi}\,
  \frac{d}{dz}
  \left[
    \frac{z}{z-z_{1}}
  \right]
  \\
  & = &
  \frac{-\ii z}{\pi}
  \left[
    \frac{1}{z-z_{1}}
    -
    \frac{z}{(z-z_{1})^{2}}
  \right]
  \\
  & = &
  \frac{-\ii z}{\pi}
  \left[
    \frac{z-z_{1}}{(z-z_{1})^{2}}
    -
    \frac{z}{(z-z_{1})^{2}}
  \right]
  \;\;\;\Rightarrow
\end{eqnarray*}
\begin{displaymath}
  w_{\delta'}(z,z_{1})
  =
  \frac{\ii zz_{1}}{\pi(z-z_{1})^{2}}.
\end{displaymath}

\noindent
It is not difficult to verify that this inner analytic function
corresponds to the sequence of Taylor-Fourier coefficients
$\alpha_{k}=k/\pi$ for all $k\geq 1$, when expanded in powers of
$z/z_{1}$, by simply taking the logarithmic derivative of
Equation~(\ref{expandelta}). Therefore, the corresponding Fourier series
also diverges almost everywhere, and does so faster than the Fourier
series of the delta ``function''. This derivation is an example of the
derivative of a singular and therefore in principle non-differentiable
object on the unit circle being regularized, and in fact defined, within
the open unit disk of the complex plane. Following the pattern of the
inner analytic function associated to the delta ``function'', which has a
simple pole on the unit circle, this one has a second-order pole at the
same point. Rationalizing this function, in order to identify its real and
imaginary parts, we get, with $z=\rho\exp(\ii\theta)$,
$z_{1}=\exp(\ii\theta_{1})$ and $\Delta\theta=\theta-\theta_{1}$,

\noindent
\begin{eqnarray*}
  \ii z\,
  \frac{d}{dz}
  w_{\delta}(z,z_{1})
  & = &
  \frac
  {
    \ii\rho
    \left(
      \rho^{2}
      \e{-\ii\Delta\theta}
      -
      2\rho
      +
      \e{\ii\Delta\theta}
    \right)
  }
  {
    \pi
    \left[
      \left(\rho^{2}+1\right)
      -
      2\rho\cos(\theta-\theta_{1})
    \right]^{2}
  }
  \\
  & = &
  \frac
  {
    \rho
    \left(\rho^{2}-1\right)
    \sin(\theta-\theta_{1})
  }
  {
    \pi
    \left[
      \left(\rho^{2}+1\right)
      -
      2\rho\cos(\theta-\theta_{1})
    \right]^{2}
  }
  +
  \\
  &   &
  +
  \ii\,
  \frac
  {
    \rho
    \left[
      -
      2\rho
      +
      \left(\rho^{2}+1\right)
      \cos(\theta-\theta_{1})
    \right]
  }
  {
    \pi
    \left[
      \left(\rho^{2}+1\right)
      -
      2\rho\cos(\theta-\theta_{1})
    \right]^{2}
  }
  \;\;\;\Rightarrow
  \\
  \delta'(\theta-\theta_{1})
  & = &
  \lim_{\rho\to 1}
  \frac
  {
    \rho
    \left(\rho^{2}-1\right)
    \sin(\theta-\theta_{1})
  }
  {
    \pi
    \left[
      \left(\rho^{2}+1\right)
      -
      2\rho\cos(\theta-\theta_{1})
    \right]^{2}
  }.
\end{eqnarray*}

\noindent
This gives us an explicit representation of the generalized function
$\delta'(\theta-\theta_{1})$, which is thus seen to depend only on
$\Delta\theta=\theta-\theta_{1}$, as a limit to the unit circle. Observe
that the expression on the right-hand side changes sign with
$\Delta\theta$, and is therefore odd with respect to $\Delta\theta$.

The process of integration by parts used to define this generalized
function can now be iterated, since we may consider the integral

\noindent
\begin{eqnarray}\label{firstform2}
  I_{\delta'}(\rho,\theta_{1})
  & = &
  \oint_{C}dz\,
  \frac{dg(\rho,\theta)}{dz}\,
  w_{\delta'}(z,z_{1})
  \nonumber\\
  & = &
  \int_{-\pi}^{\pi}d\theta\,
  \frac{dg(\rho,\theta)}{d\theta}\,
  w_{\delta'}(z,z_{1})
  \nonumber\\
  & = &
  \int_{-\pi}^{\pi}d\theta\,
  \frac{dg(\rho,\theta)}{d\theta}
  \left[
    f'(\rho,\theta,\theta_{1})
    +
    \ii
    \bar{f}'(\rho,\theta,\theta_{1})
  \right],
\end{eqnarray}

\noindent
which once more can be integrated by parts to yield

\noindent
\begin{eqnarray}\label{secondform2}
  I_{\delta'}(\rho,\theta_{1})
  & = &
  -
  \oint_{C}dz\,
  g(\rho,\theta)\,
  \frac{d}{dz}
  w_{\delta'}(z,z_{1})
  \nonumber\\
  & = &
  -
  \int_{-\pi}^{\pi}d\theta\,
  g(\rho,\theta)\,
  \left[
    \frac{d}{d\theta}
    f'(\rho,\theta,\theta_{1})
    +
    \ii\,
    \frac{d}{d\theta}
    \bar{f}'(\rho,\theta,\theta_{1})
  \right].
\end{eqnarray}

\noindent
We now observe that the real part of the original form of this integral in
Equation~(\ref{firstform2}), when the limit to the unit circle is taken,
becomes the integral of the product of the first derivative of the delta
``function'' and the derivative of the real function, so that we have

\noindent
\begin{eqnarray*}
  I_{\delta'}(1,\theta_{1})
  & = &
  \int_{-\pi}^{\pi}d\theta\,
  \left[
    \frac{dg(1,\theta)}{d\theta}\,
    \delta'(\theta-\theta_{1})
    +
    \ii\,
    \frac{dg(1,\theta)}{d\theta}\,
    \bar{f}'(1,\theta,\theta_{1})
  \right]
  \\
  & = &
  -\,
  \frac{d^{2}g}{d\theta^{2}}(1,\theta_{1})
  +
  \ii
  \int_{-\pi}^{\pi}d\theta\,
  \frac{dg(1,\theta)}{d\theta}\,
  \bar{f}'(1,\theta,\theta_{1}),
\end{eqnarray*}

\noindent
where we used the integration properties of $\delta'(\theta-\theta_{1})$.
If we also take the limit to the unit circle of the last expression in
Equation~(\ref{secondform2}), then it follows that on that circle we have

\noindent
\begin{eqnarray*}
  I_{\delta'}(1,\theta_{1})
  & = &
  -
  \lim_{\rho\to 1}
  \int_{-\pi}^{\pi}d\theta\,
  g(\rho,\theta)\,
  \left[
    \frac{d}{d\theta}
    f'(\rho,\theta,\theta_{1})
    +
    \ii\,
    \frac{d}{d\theta}
    \bar{f}'(\rho,\theta,\theta_{1})
  \right]
  \\
  & = &
  -
  \int_{-\pi}^{\pi}d\theta\,
  g(1,\theta)
  \lim_{\rho\to 1}
  \left[
    \frac{d}{d\theta}
    f'(\rho,\theta,\theta_{1})
  \right]
  -
  \ii
  \int_{-\pi}^{\pi}d\theta\,
  g(1,\theta)\,
  \lim_{\rho\to 1}
  \left[
    \frac{d}{d\theta}
    \bar{f}'(\rho,\theta,\theta_{1})
  \right].
\end{eqnarray*}

\noindent
From these two expressions for $I_{\delta'}(1,\theta_{1})$ we conclude
that we have

\noindent
\begin{eqnarray*}
  \lefteqn
  {
    -\,
    \frac{d^{2}g}{d\theta^{2}}(1,\theta_{1})
    +
    \ii
    \int_{-\pi}^{\pi}d\theta\,
    \frac{dg(1,\theta)}{d\theta}\,
    \bar{f}'(1,\theta,\theta_{1})
  }
  &   &
  \\
  & = &
  -
  \int_{-\pi}^{\pi}d\theta\,
  g(1,\theta)
  \lim_{\rho\to 1}
  \left[
    \frac{d}{d\theta}
    f'(\rho,\theta,\theta_{1})
  \right]
  -
  \ii
  \int_{-\pi}^{\pi}d\theta\,
  g(1,\theta)\,
  \lim_{\rho\to 1}
  \left[
    \frac{d}{d\theta}
    \bar{f}'(\rho,\theta,\theta_{1})
  \right].
\end{eqnarray*}

\noindent
This time the real part of this equality tell us that the expression
involving the limit of the derivative shown within the integral at the
right-hand side can be interpreted as a definition of the second
derivative of the delta ``function'',

\begin{displaymath}
  \frac{d^{2}g}{d\theta^{2}}(1,\theta_{1})
  =
  \int_{-\pi}^{\pi}d\theta\,
  g(1,\theta)
  \lim_{\rho\to 1}
  \left[
    \frac{d}{d\theta}
    f'(\rho,\theta,\theta_{1})
  \right]
  \;\;\;\Rightarrow
\end{displaymath}
\begin{displaymath}
  \frac{d^{2}}{d\theta^{2}}
  \delta(\theta-\theta_{1})
  \equiv
  \lim_{\rho\to 1}
  \left[
    \frac{d}{d\theta}
    f'(\rho,\theta,\theta_{1})
  \right].
\end{displaymath}

\noindent
This is a new generalized function, that has the property that through
integration over the unit circle it attributes to a given
twice-differentiable function the value of its second derivative at the
point $\theta_{1}$,

\begin{displaymath}
  \int_{-\pi}^{\pi}d\theta\,
  g(1,\theta)\,
  \frac{d^{2}}{d\theta^{2}}
  \delta(\theta-\theta_{1})
  =
  \frac{d^{2}g}{d\theta^{2}}(1,\theta_{1}).
\end{displaymath}

\noindent
Once again we see that this is the integration kernel of a real-valued
linear functional acting in the space of integrable real functions. We
will adopt for this new singular object the notation

\begin{displaymath}
  \delta''(\theta-\theta_{1})
  \equiv
  \frac{d^{2}}{d\theta^{2}}
  \delta(\theta-\theta_{1}).
\end{displaymath}

\noindent
We may now determine the inner analytic function associated to this new
singular object. Once again, since derivatives with respect to $\theta$
may be written as logarithmic derivatives of the corresponding inner
analytic function within the open unit disk, as shown in~\cite{FTotCPII},
we have, using the notation established in that paper,

\noindent
\begin{eqnarray*}
  w_{\delta''}(z,z_{1})
  & = &
  w^{\ldot}_{\delta'}(z,z_{1})
  \\
  & = &
  \ii z\,
  \frac{d}{dz}
  w_{\delta'}(z,z_{1})
  \\
  & = &
  \frac{d}{d\theta}
  f'(\rho,\theta,\theta_{1})
  +
  \ii
  \frac{d}{d\theta}
  \bar{f}'(\rho,\theta,\theta_{1}),
\end{eqnarray*}

\noindent
so that the generalized function $\delta''(\theta-\theta_{1})$ is given by
the limit to the unit circle of the inner analytic function
$w_{\delta''}(z,z_{1})$,

\begin{displaymath}
  \delta''(\theta-\theta_{1})
  =
  \lim_{\rho\to 1}
  \Re[w_{\delta''}(z,z_{1})].
\end{displaymath}

\noindent
Since we have the inner analytic function that corresponds to the first
derivative of the delta ``function'' in explicit form, we have the
capability to perform an explicit calculation in order to get the inner
analytic function that corresponds to the second derivative of the delta
``function''. It can be easily verified that the inner analytic function
$w_{\delta''}(z,z_{1})$ corresponds to the sequence of Taylor-Fourier
coefficients $\alpha_{k}=k^{2}/\pi$ for all $k\geq 1$, when expanded in
powers of $z/z_{1}$.

Just as the first and the second logarithmic derivatives of the inner
analytic function corresponding to the original delta ``function''
represent the first and second derivatives of the delta ``function''
within the open unit disk, so it is with all the higher derivatives. By
iterating repeatedly this process of integration by parts, we can generate
a whole infinite sequence of singular generalized functions and their
corresponding inner analytic functions. Iterating one more time one gets
for the first few elements of this sequence of inner analytic functions
the list of results that follows,

\noindent
\begin{eqnarray*}
  w_{\delta^{1\prime}}(z,z_{1})
  & = &
  -\,
  \frac{1}{\pi\ii^{1}}\,
  \frac
  {zz_{1}}
  {(z-z_{1})^{2}},
  \\
  w_{\delta^{2\prime}}(z,z_{1})
  & = &
  -\,
  \frac{1}{\pi\ii^{2}}\,
  \frac
  {z(z+z_{1})z_{1}}
  {(z-z_{1})^{3}},
  \\
  w_{\delta^{3\prime}}(z,z_{1})
  & = &
  -\,
  \frac{1}{\pi\ii^{3}}\,
  \frac
  {z\left(z^{2}+4zz_{1}+z_{1}^{2}\right)z_{1}}
  {(z-z_{1})^{4}}.
\end{eqnarray*}

\noindent
Generalizing to the $n^{\rm th}$ derivative, we see that the manifestly
systematic factors in $w_{\delta^{n\prime}}(z,z_{1})$ can be written as

\begin{displaymath}
  -\,
  \frac{zz_{1}}{\pi\ii^{n}(z-z_{1})^{n+1}},
\end{displaymath}

\noindent
and the remaining factors are polynomials on $z$. Note that the
logarithmic derivation operator conserves the homogeneity of the powers of
$z$ and $z_{1}$, because we have that

\begin{displaymath}
  \ii z\,\frac{d}{dz}\blank
  =
  \ii(z/z_{1})\,\frac{d}{d(z/z_{1})}\blank.
\end{displaymath}

\noindent
Therefore all these inner analytic functions can be written as functions
of only $z/z_{1}$, so that they are all rotated inner analytic functions,
with the angle $\theta_{1}$. Also, it is not difficult to verify that
$w_{\delta^{n\prime}}(z,z_{1})$ corresponds to the sequence of
Taylor-Fourier coefficients $\alpha_{k}=k^{n}/\pi$ for all $k\geq 1$, when
expanded in powers of $z/z_{1}$. It is possible to work out a general
recursion formula for the inner analytic functions corresponding to all
the multiple derivatives of the delta ``function''. In order to do this we
may write, in terms of the variable $\chi=z/z_{1}$, for $n>0$, that

\begin{displaymath}
  w_{\delta^{n\prime}}(\chi)
  =
  -\,
  \frac{1}{\pi\ii^{n}}\,
  \frac{\chi P_{n-1}(\chi)}{(\chi-1)^{n+1}},
\end{displaymath}

\noindent
where $P_{n-1}(\chi)$ is a polynomial of order $n-1$ on $\chi$, with real
(in fact, integer) coefficients $A_{n,i}$,

\begin{displaymath}
  P_{n}(\chi)
  =
  \sum_{i=0}^{n}
  A_{n,i}\chi^{i}.
\end{displaymath}

\noindent
One can see that the first few polynomials, for $n>0$, are therefore given
by

\noindent
\begin{eqnarray*}
  P_{0}(\chi)
  & = &
  1,
  \\
  P_{1}(\chi)
  & = &
  1+\chi,
  \\
  P_{2}(\chi)
  & = &
  1+4\chi+\chi^{2}.
\end{eqnarray*}

\noindent
As is shown in Appendix~\ref{APPcalcoefs}, one obtains for the
coefficients $A_{n,i}$ the purely algebraic two-index recurrence relation

\begin{table}
\noindent
\begin{displaymath}
  \begin{array}{|l|rrrrrrr|}
    \hline
    i=\to  & \ \  0 & \ \  1 & \ \  2 & \ \  3 & \ \  4 & \ \  5 & \cdots \\
    \hline
    n=0    &      1 &      0 &      0 &      0 &      0 &      0 & \cdots \\
    n=1    &      1 &      1 &      0 &      0 &      0 &      0 & \cdots \\
    n=2    &      1 &      4 &      1 &      0 &      0 &      0 & \cdots \\
    n=3    &      1 &     11 &     11 &      1 &      0 &      0 & \cdots \\
    n=4    &      1 &     26 &     66 &     26 &      1 &      0 & \cdots \\
    n=5    &      1 &     57 &    302 &    302 &     57 &      1 & \cdots \\
    \vdots & \vdots & \vdots & \vdots & \vdots & \vdots & \vdots & \ddots \\
    \hline
  \end{array}
\end{displaymath}
\caption{A table of values for the coefficients $A_{n,i}$.}\label{thetable}
\end{table}

\begin{equation}\label{algebrecur}
  A_{n,i}
  =
  (i+1)A_{n-1,i}+(n-i+1)A_{n-1,i-1},
\end{equation}

\noindent
for $0<i<n$. When using this formula we should assume that $A_{n,i}=0$ for
$i<0$ and for $i>n$, and we must also use the initial value $A_{0,0}=1$.
Using this recursion relation we may construct the table of values of
$A_{n,i}$ shown in Table~\ref{thetable}. Using this method we are
therefore able to determine completely the inner analytic function
associated to the derivative of any given order of the delta ``function''.

Note that the $n^{\rm th}$ derivative of the delta ``function'' is
represented by an inner analytic function with a pole of order $n+1$ on
the unit circle, and that, as we discussed before, it is associated to the
sequence of Taylor-Fourier coefficients $\alpha_{k}=k^{n}/\pi$ for all
$k\geq 1$, when expanded in powers of $\chi$. We may therefore conclude
that all these generalized functions are within the set defined by the
limited definition discussed in the introduction, and expressed by
Equations~(\ref{limitcond1}) and~(\ref{limitcond2}), including therefore
the requirement that the sequence of Taylor-Fourier coefficients
$\alpha_{k}$ not diverge faster than all powers of $k$ when $k\to\infty$.

\section{A Class of Locally Non-Integrable Functions}\label{SEClocalnon}

In we consider once again the extended inner analytic function
$w_{\delta}(z,z_{1})$ associated to the delta ``function'', given in
Equation~(\ref{innerdelta}), we realize that not only its real part is a
representation of the delta ``function'' in the $\rho\to 1$ limit, but the
corresponding FC function is also representable by the same extended inner
analytic function, as the $\rho\to 1$ limit of its imaginary part. This is
a very simple normal function $\bar{f}(\theta-\theta_{1})$, although a
singular and non-integrable one, which was given in the appendices
of~\cite{FTotCPI}, and which we repeat here,

\begin{displaymath}
  \bar{f}(\theta-\theta_{1})
  =
  \frac{1}{\pi}\,
  \frac{1+\cos(\theta-\theta_{1})}{2\sin(\theta-\theta_{1})}.
\end{displaymath}

\noindent
This function has a singularity that behaves as $1/(\theta-\theta_{1})$
around the point $\theta_{1}$, which is not, therefore, an integrable
singularity. This means that the integral of this function over the unit
circle is not well defined. However, it is still associated to the
Taylor-Fourier coefficients of the delta ``function'', and can be
recovered as the $\rho\to 1$ limit of the imaginary part of the same
extended inner analytic function. We may also recover this function as the
real part of its own extended inner analytic function, which is simply
given by

\begin{displaymath}
  w(z,z_{1})
  =
  -\ii
  w_{\delta}(z,z_{1}),
\end{displaymath}

\noindent
being thus defined according to our new criterion. It follows from this
fact that it is not really necessary for a real function to be integrable
in order for it to have a well-defined sequence of Taylor-Fourier
coefficients and to be recoverable almost everywhere from, and therefore
representable almost everywhere by, an inner analytic function. In fact,
it is possible to identify without too much trouble a whole class of
singular non-integrable real functions that can be represented by inner
analytic functions. All we have to do in order to accomplish this is to
find a consistent way to associate a sequence of Fourier coefficients to
such functions.

We will call functions that are locally integrable almost everywhere, by
which we mean everywhere but in the neighborhoods of a finite set of
singular points on the unit circle, by the name of {\em locally
  non-integrable functions}. This means that the function is integrable in
all closed intervals within the unit circle that do {\em not} contain a
singular point, of which we assume there is a finite number. We will be
able to define Fourier coefficients for such functions so long as they can
be obtained as derivatives of any finite order of integrable functions.
Let us start the argument with the case of the first derivative. If a real
function $f(\theta)$ is locally non-integrable, so that the coefficients

\noindent
\begin{eqnarray*}
  \alpha_{k}
  & = &
  \frac{1}{\pi}
  \int_{-\pi}^{\pi}d\theta\,
  \cos(k\theta)
  f(\theta),
  \\
  \beta_{k}
  & = &
  \frac{1}{\pi}
  \int_{-\pi}^{\pi}d\theta\,
  \sin(k\theta)
  f(\theta),
\end{eqnarray*}

\noindent
for $k\geq 1$, do not exist if defined in this way, due to the fact that
the integrals do not exist, but $f(\theta)$ is the derivative with respect
to $\theta$ of a zero-average real function $f^{-1\prime}(\theta)$ which
is defined almost everywhere and integrable on the whole unit circle, then
we may define the Taylor-Fourier coefficients associated to $f(\theta)$ as

\noindent
\begin{eqnarray*}
  \alpha_{k}
  & = &
  \frac{k}{\pi}
  \int_{-\pi}^{\pi}d\theta\,
  \sin(k\theta)
  f^{-1\prime}(\theta),
  \\
  \beta_{k}
  & = &
  -\,
  \frac{k}{\pi}
  \int_{-\pi}^{\pi}d\theta\,
  \cos(k\theta)
  f^{-1\prime}(\theta),
\end{eqnarray*}

\noindent
for $k\geq 1$, where we note that the sine and cosine have been
interchanged. We might call this the {\em extended definition} of the
Fourier coefficients. Note that if $f(\theta)$ is the derivative of a
zero-average function, then it also has no constant term in its Fourier
expansion, and thus we may consider it to be zero-average as well. Note
also that when $f(\theta)$ is integrable and the usual definitions of
$\alpha_{k}$ and $\beta_{k}$ are in effect, we may always write that

\noindent
\begin{eqnarray*}
  \alpha_{k}
  & = &
  \frac{1}{\pi}
  \int_{-\pi}^{\pi}d\theta\,
  \cos(k\theta)
  f(\theta)
  \\
  & = &
  \frac{k}{\pi}
  \int_{-\pi}^{\pi}d\theta\,
  \sin(k\theta)
  f^{-1\prime}(\theta),
  \\
  \beta_{k}
  & = &
  \frac{1}{\pi}
  \int_{-\pi}^{\pi}d\theta\,
  \sin(k\theta)
  f(\theta)
  \\
  & = &
  -\,
  \frac{k}{\pi}
  \int_{-\pi}^{\pi}d\theta\,
  \cos(k\theta)
  f^{-1\prime}(\theta),
\end{eqnarray*}

\noindent
for $k\geq 1$, since we can always integrate by parts on the unit circle,
without generating an integrated term. The point is that even if
$f(\theta)$ is not integrable and the usual definition does not work, the
extended definition may still exist. This will be the case when
$f(\theta)$, though not an integrable function itself, is the derivative
of an integrable function. This will often happen around singular points
where $f(\theta)$ diverges to infinity. Another way to implement this
extended definition is to write the Fourier coefficients $\alpha_{k}$ and
$\beta_{k}$ of the function $f(\theta)$ in terms of the Fourier
coefficients $\alpha^{-1\prime}_{k}$ and $\beta^{-1\prime}_{k}$ of its
zero-average first primitive $f^{-1\prime}(\theta)$,

\noindent
\begin{eqnarray*}
  \alpha_{k}
  & = &
  k
  \beta_{k}^{-1\prime},
  \\
  \beta_{k}
  & = &
  -k
  \alpha_{k}^{-1\prime}.
\end{eqnarray*}

\noindent
The best way to interpret this extended definition is in terms of the
corresponding inner analytic functions within the open unit disk. In order
to do this, one starts with a locally non-integrable function $f(\theta)$.
Since it is integrable in all closed intervals within the unit circle that
do {\em not} contain a singularity, one may in principle integrate within
these intervals in order to find a zero-average primitive
$f^{-1\prime}(\theta)$. It should be possible, in this way, to determine a
zero-average real function $f^{-1\prime}(\theta)$ such that $f(\theta)$ is
its derivative almost everywhere. So long as $f^{-1\prime}(\theta)$ is
integrable, one may then calculate its Taylor-Fourier coefficients in the
usual way, and thus construct the corresponding inner analytic function
$w^{-1\!\ldot}(z)$. By taking the logarithmic derivative of this function
one then gets the inner analytic function $w(z)$ that corresponds to
$f(\theta)$, and at the same time determines its Taylor-Fourier
coefficients.

Note that if $f(\theta)$ is not integrable then the Fourier series
associated to its extended Fourier coefficients will typically diverge
almost everywhere, due to the extra factor of $k$ in the coefficients.
However, so long as the coefficients do not diverge faster than all powers
of $k$, the corresponding inner analytic function can still be constructed
with those coefficients, and the function can still be recovered from it
almost everywhere. Therefore, it is still true that the locally
non-integrable real function is uniquely characterized by its extended
Fourier coefficients almost everywhere.

The definition of the Fourier coefficients can be further extended to
cases in which both $f(\theta)$ and $f^{-1\prime}(\theta)$ are locally
non-integrable, so long as $f(\theta)$ is the second derivative with
respect to $\theta$ of a zero-average integrable real function
$f^{-2\prime}(\theta)$. In this case we may define

\noindent
\begin{eqnarray*}
  \alpha_{k}
  & = &
  -\,
  \frac{k^{2}}{\pi}
  \int_{-\pi}^{\pi}d\theta\,
  \sin(k\theta)
  f^{-2\prime}(\theta),
  \\
  \beta_{k}
  & = &
  -\,
  \frac{k^{2}}{\pi}
  \int_{-\pi}^{\pi}d\theta\,
  \cos(k\theta)
  f^{-2\prime}(\theta),
\end{eqnarray*}

\noindent
for $k\geq 1$, by simply using for $f^{-1\prime}(\theta)$ the same
argument used previously for $f(\theta)$. Note that the sine and cosine
have now been interchanged back to their original positions. With the use
of repeated integrations by parts this can be extended to locally
non-integrable functions $f(\theta)$ which are the derivatives of any
finite order of integrable real functions with finite numbers of
singularities where they diverge to infinity. It is not difficult to
derive general formulas for the coefficients $\alpha_{k}$ and $\beta_{k}$
valid for the case of the derivative of order $n$, and the corresponding
partial integrations until the primitive of order $n$ is produced. We just
have to systematize the first few cases, which are given by

\noindent
\begin{eqnarray*}
  \alpha_{k}
  & = &
  +\,
  \frac{k^{0}}{\pi}
  \int_{-\pi}^{\pi}d\theta\,
  \cos(k\theta)
  f^{-0\prime}(\theta)
  \\
  & = &
  +\,
  \frac{k^{1}}{\pi}
  \int_{-\pi}^{\pi}d\theta\,
  \sin(k\theta)
  f^{-1\prime}(\theta)
  \\
  & = &
  -\,
  \frac{k^{2}}{\pi}
  \int_{-\pi}^{\pi}d\theta\,
  \cos(k\theta)
  f^{-2\prime}(\theta)
  \\
  & = &
  -\,
  \frac{k^{3}}{\pi}
  \int_{-\pi}^{\pi}d\theta\,
  \sin(k\theta)
  f^{-3\prime}(\theta)
  \\
  & = &
  +\,
  \frac{k^{4}}{\pi}
  \int_{-\pi}^{\pi}d\theta\,
  \cos(k\theta)
  f^{-4\prime}(\theta)
  \\
  & = &
  +\,
  \frac{k^{5}}{\pi}
  \int_{-\pi}^{\pi}d\theta\,
  \sin(k\theta)
  f^{-5\prime}(\theta),
  \\
  \beta_{k}
  & = &
  +\,
  \frac{k^{0}}{\pi}
  \int_{-\pi}^{\pi}d\theta\,
  \sin(k\theta)
  f^{-0\prime}(\theta)
  \\
  & = &
  -\,
  \frac{k^{1}}{\pi}
  \int_{-\pi}^{\pi}d\theta\,
  \cos(k\theta)
  f^{-1\prime}(\theta)
  \\
  & = &
  -\,
  \frac{k^{2}}{\pi}
  \int_{-\pi}^{\pi}d\theta\,
  \sin(k\theta)
  f^{-2\prime}(\theta)
  \\
  & = &
  +\,
  \frac{k^{3}}{\pi}
  \int_{-\pi}^{\pi}d\theta\,
  \cos(k\theta)
  f^{-3\prime}(\theta)
  \\
  & = &
  +\,
  \frac{k^{4}}{\pi}
  \int_{-\pi}^{\pi}d\theta\,
  \sin(k\theta)
  f^{-4\prime}(\theta)
  \\
  & = &
  -\,
  \frac{k^{5}}{\pi}
  \int_{-\pi}^{\pi}d\theta\,
  \cos(k\theta)
  f^{-5\prime}(\theta),
\end{eqnarray*}

\noindent
where we used the notation $f(\theta)=f^{-0\prime}(\theta)$. In order to
systematize these cases we must separate then into those in which $n$ is
even and those in which $n$ is odd. For even $n$ we make $n=2j$, with
$j=0,1,2,\ldots,\infty$, and we may write

\noindent
\begin{eqnarray*}
  \alpha_{k}
  & = &
  \frac{(-1)^{j}k^{n}}{\pi}
  \int_{-\pi}^{\pi}d\theta\,
  \cos(k\theta)
  f^{-n\prime}(\theta),
  \\
  \beta_{k}
  & = &
  \frac{(-1)^{j}k^{n}}{\pi}
  \int_{-\pi}^{\pi}d\theta\,
  \sin(k\theta)
  f^{-n\prime}(\theta),
\end{eqnarray*}

\noindent
for $k\geq 1$, while for odd $n$ we make $n=2j+1$, with
$j=0,1,2,\ldots,\infty$, and we may write

\noindent
\begin{eqnarray}\label{extendefcoef}
  \alpha_{k}
  & = &
  \frac{(-1)^{j}k^{n}}{\pi}
  \int_{-\pi}^{\pi}d\theta\,
  \sin(k\theta)
  f^{-n\prime}(\theta),
  \nonumber\\
  \beta_{k}
  & = &
  \frac{(-1)^{j+1}k^{n}}{\pi}
  \int_{-\pi}^{\pi}d\theta\,
  \cos(k\theta)
  f^{-n\prime}(\theta),
\end{eqnarray}

\noindent
for $k\geq 1$, where we note the interchange of the sine and the cosine,
as well as the extra sign in the second equation. These equations
constitute the most general definition of the {\em extended Taylor-Fourier
  coefficients} of the function $f(\theta)$. The idea is that we use for
the calculation of the coefficients the smallest value of $n$ for which
the integrals exist, that is, the smallest value of $n$ for which the
primitive of order $n$ of $f(\theta)$, that is the function
$f^{-n\prime}(\theta)$, is a zero-average integrable function.

Since the Fourier coefficients of an integrable function are always
limited, it is clear that the extended Fourier coefficients of the class
of locally non-integrable functions that we are considering here, which
are given by those limited coefficients multiplied by a fixed power of
$k$, will never diverge to infinity faster than all powers of $k$ as we
make $k\to\infty$. This is in accordance of our restrictive hypothesis,
given in Equation~(\ref{limitcond2}), and hence we see that with the use
of the extended coefficients we are always able to represent this class of
locally non-integrable functions by means of inner analytic functions
within the open unit disk. We are, therefore, always able to recover these
functions almost everywhere as limits to the unit circle of the real part
of the corresponding inner analytic functions, although the corresponding
Fourier series are always divergent almost everywhere.

As a final note, let us observe that, although the extended Fourier
coefficients of non-integrable functions will always result in divergent
Fourier series, using the analytic structure within the open unit disk one
may apply to these coefficients the technique of {\em singularity
  factorization} developed in~\cite{FTotCPII}, in order to obtain
convergent trigonometric series to represent these functions everywhere on
the unit circle but at the singular points. In fact, several such series
can be constructed from the extended Fourier coefficients, with increasing
levels of speed of convergence. These are the {\em center series} which
were developed and described in~\cite{FTotCPII}, and tested numerically
in~\cite{CSNumerics}. In this way, the extended Fourier coefficients are
seen to still have an algorithmic value, in terms of the numerical
representation of locally non-integrable functions, despite the fact that
the corresponding Fourier series are divergent.

\section{Boundary Value Problems on the Unit Disk}\label{SECboundval}

Let us show that the correspondence between inner analytic functions
within the open unit disk and real functions on the unit circle can be
interpreted in terms of boundary value problems of the two-dimensional
Laplace equation on the unit disk. In order to do this, let us describe
how one goes about finding the inner analytic function that corresponds to
a given integrable real function. We will formulate the ideas in terms of
normal integrable real functions, and later generalize the results that
are found. The idea is to develop a criterion for the existence of the
inner analytic function associated to a normal integrable real function,
that can later be applied to both normal and generalized function, as well
as to both regular and singular functions.

Let an arbitrary integrable real function $f(\theta)$ be given. According
to the results we obtained before, the corresponding inner analytic
function can be constructed in the following way: since this real function
is integrable on the unit circle, one can calculate the corresponding
Taylor-Fourier coefficients and then use them to build the corresponding
inner analytic function by means of a convergent complex Taylor series.
This produces an inner analytic function that can be written as

\begin{displaymath}
  w(z)
  =
  f(\rho,\theta)
  +
  \ii
  \bar{f}(\rho,\theta),
\end{displaymath}

\noindent
where $f(\rho,\theta)$ is harmonic and $\bar{f}(\rho,\theta)$ is its
harmonic conjugate. The process described above can be understood as a way
to determine $f(\rho,\theta)$ from $f(\theta)$. If $f(\rho,\theta)$ can be
defined as a harmonic function on the open unit disk, then its harmonic
conjugate $\bar{f}(\rho,\theta)$ necessarily exists and is harmonic on the
same domain. This is so because a harmonic function defined on a simply
connected domain within the complex plane always admits a harmonic
conjugate. In fact, in our case here the harmonic conjugate is given by a
real line integral within the open unit disk, written here both in
Cartesian coordinates and in a somewhat more abstract way using vector
notation,

\noindent
\begin{eqnarray*}
  \bar{f}(x,y)
  & = &
  \int_{(0,0)}^{(x,y)}
  \left[
    \frac{\partial f\!\left(x',y'\right)}{\partial x'}\,dy'
    -
    \frac{\partial f\!\left(x',y'\right)}{\partial y'}\,dx'
  \right]
  \\
  & = &
  \int_{0}^{z}
  \left[
    \vec{\nabla}f\!\left(z'\right)\times\vec{dz'}
  \right]_{n},
\end{eqnarray*}

\noindent
where what we have here is actually the component of a vector product
normal to the complex plane, and where the infinitesimal vector is defined
as $\vec{dz}=(dx,dy,0)$. Therefore, we may focus our attention on the
determination of $f(\rho,\theta)$, aiming at describing it in a more
general way, so as not to depend on the determination of the Fourier
coefficients by means of integrals over the unit circle. Now, the problem
of finding $f(\rho,\theta)$ within the open unit disk starting from the
values of $f(\theta)$ on the unit circle can be understood as a standard
{\em boundary value problem} on the unit disk. In fact, since the real and
imaginary parts of an analytic function are harmonic functions, it is the
problem of finding a solution of the two-dimensional Laplace equation,
written here in polar coordinates,

\begin{displaymath}
  \frac{\partial^{2}f(\rho,\theta)}{\partial\rho^{2}}
  +
  \frac{1}{\rho^{2}}\,
  \frac{\partial^{2}f(\rho,\theta)}{\partial\theta^{2}}
  =
  0,
\end{displaymath}

\noindent
valid within the open unit disk, that satisfies the boundary condition,
over the whole unit circle,

\begin{displaymath}
  f(1,\theta)
  =
  f(\theta).
\end{displaymath}

\noindent
This is an instance of the famous Dirichlet problem, to wit the Dirichlet
problem on the unit disk~\cite{Dirichlet}. The simplest theorem of
existence and uniqueness of the solution of boundary value problems such
as this one states that the solution exists and is unique so long as the
boundary condition $f(\theta)$ is well defined and continuous over the
whole boundary. There are further theorems that establish the existence
and uniqueness of the solution for boundary conditions $f(\theta)$ that
are well defined everywhere but discontinuous at some points. There seems
to be no well-established general theorems about the solution of the
boundary value problem when the real function $f(\theta)$ has points of
singularity where it diverges to infinity.

The existence or non-existence of the solution of the boundary value
problem is a very general way to determine whether or not the inner
analytic function can be defined. We may consider applying it not only to
integrable normal real functions, but also to locally non-integrable
normal real functions, as well as to singular generalized functions. In
these two singular cases the main difference is the existence of points of
singularity where the quantity $f(\theta)$ diverges to infinity. In these
cases it is not enough to apply as boundary conditions the known values of
the functions at the unit circle, since this clearly cannot be done at the
singular points. The set of boundary conditions are thus seen to be
incomplete in such cases. If we are to have a complete set of boundary
conditions, the boundary condition at every point of divergence of
$f(\theta)$ must be replaced by an individual extra condition, for example
a condition such as the one about the integral of the delta ``function''.

This leads one to think of a modified boundary value problem, in which the
continuously infinite set of conditions on the boundary over the whole
unit circle is replaced by a boundary condition with a certain number of
point-like ``holes''. We might call this new type of boundary conditions
``punctured boundary conditions''. In such cases one must add to the
boundary-value problem a set of extra conditions, in the same number as
the holes. In the case of normal real functions that diverge to infinity
at the singular points, the extra conditions might be statements about the
form or speed of the divergence. This implies that the behavior of the
real functions in neighborhoods around the singular points is known and
can be used in this way, as local conditions around each singular point.

For more singular objects such as the delta ``function'' the value of the
``function'' in such neighborhoods is simply zero, except at the singular
points themselves, so that information regarding how the ``function''
tends to infinity is not available on the unit circle and thus cannot be
used to characterize the ``functions''. Therefore, one is left only with
the alternative of imposing a global rather than local extra condition,
such as the one involving the integral of the delta ``function''. Seen in
this light, the Dirac delta ``function'', and the other singular
generalized functions associated to it that we examined in
Section~\ref{SECgenfunc}, are all solutions of punctured boundary value
problems of the Laplace equation on the unit disk, with a single singular
point $z_{1}$ on the unit circle.

We are led therefore to the fact that there is an interesting side-effect
of the correspondence between inner analytic functions within the open
unit disk and real objects on the unit circle. The structure we introduced
at once implies some quite general existence and uniqueness theorems for
the solution of the Dirichlet problem on the unit disk. For example, from
the fact that every integrable real function has a unique corresponding
inner analytic function we may conclude that the solution of the Dirichlet
problem on the unit disk exists and is unique for any boundary condition
$f(\theta)$ that is simply integrable on the unit circle. This extends the
validity of the usual theorem to functions $f(\theta)$ that not only can
be discontinuous, but that can also have integrable divergences at a
finite number of singular points.

Since we have extended the set of inner analytic functions within the open
unit disk to include on the unit circle a class of locally non-integrable
real functions, the validity of the theorem is also extended to that class
of singular boundary conditions, with a suitable reinterpretation of the
set of boundary conditions. Finally, since we have extended the set of
inner analytic functions within the open unit disk to include on the unit
circle a class of singular generalized functions, the validity of the
theorem is extended to that class of singular boundary conditions as well,
with another suitable reinterpretation of the set of boundary conditions.
Given the involvement of the Dirac delta ``function'', this is clearly
related to the Green's-function approach to the solution of partial
differential equations.

\section{The Role of Integral-Differential Chains}

In~\cite{FTotCPII} we introduced the concept of an integral-differential
chain of functions. This was done for the analysis of the convergence
issues of the associated series on the unit circle. Associated to this, in
that same paper we introduced the concepts of soft and hard singularities,
as well as a corresponding gradation, in the form of the definition of a
degree of softness or of a degree of hardness of any given soft or hard
singularity. Let us now review the concept of an integral-differential
chain in the context of our new criterion for defining inner analytic
functions.

As we discussed in Section~\ref{SECgenreal}, given any integrable
zero-average real function $f(\theta)$, we can construct from it an inner
analytic function $w(z)$ such that the real function is recovered from its
real part in the $\rho\to 1$ limit. Now, once we have this complex
function, since it is analytic within the open unit disk, and hence both
infinitely differentiable and infinitely integrable there, we can consider
its successive logarithmic derivatives and logarithmic primitives, all of
which exist and are equally analytic, in exactly the same domain. Both
these operations were defined in~\cite{FTotCPII}, the logarithmic
derivative as

\begin{displaymath}
  w^{\ldot}(z)
  =
  z\,
  \frac{dw(z)}{dz},
\end{displaymath}

\noindent
which also gives our notation for it, and the logarithmic primitive as

\begin{displaymath}
  w^{-1\!\ldot}(z)
  =
  \int_{0}^{z}dz'\,
  \frac{1}{z'}\,
  w(z'),
\end{displaymath}

\noindent
over any simple integration contour contained within the open unit disk,
going from $0$ to $z$. These two operations are the inverses of one
another, as shown in~\cite{FTotCPII}. Note that the operation of
logarithmic integration was defined in a way that eliminates the usual
arbitrary integration constant, and that both definitions keep invariant
the property that $w(0)=0$. In other words, both the logarithmic
derivative and the logarithmic primitive of an inner analytic function are
also inner analytic functions, according to our new definition.

We are thus able to define a discrete infinite chain of inner analytic
functions, passing through the inner analytic function we started with,
and running by logarithmic differentiation in one direction and by
logarithmic integration in the other, indefinitely in both directions.
Since the operations of logarithmic differentiation and logarithmic
integration always produce definite and unique results, every zero-average
integrable real function belongs to a single one of these chains. Other
ways to state this are to say that two different chains cannot have a
common element, or that two chains cannot ``cross''. They constitute a
type of ``discrete fibration'' of the space of inner analytic functions.
Also, since the operations of logarithmic differentiation and logarithmic
integration correspond respectively to differentiation and integration
with respect to $\theta$ on the unit circle, given that we have on any
circle with radius $\rho\leq 1$ centered at the origin, as shown in
in~\cite{FTotCPII},

\noindent
\begin{eqnarray*}
  \frac{dw(z)}{d\theta}
  & = &
  \ii w^{\ldot}(z),
  \\
  \int_{z_{0}}^{z}d\theta'\,
  w(z')
  & = &
  -\ii
  \left[
    w^{-1\!\ldot}(z)
    -
    w^{-1\!\ldot}(z_{0})
  \right],
\end{eqnarray*}

\noindent
it follows that in the $\rho\to 1$ limit the real part of each one of the
inner analytic functions within a chain corresponds either to a derivative
$f^{n\prime}(\theta)$ or to a primitive $f^{-n\prime}(\theta)$ of the
original zero-average real function $f(\theta)$. We conclude therefore
that, given an integral-differential chain of inner analytic functions, a
corresponding chain of real objects on the unit circle is also defined. As
we discussed in Section~\ref{SECgenfunc}, these objects may or may not be
integrable zero-average real functions. Whatever the case may be, we will
from now on consider these objects to be an integral part of the chain,
forming its real sector, as illustrated by the diagram that follows:

\begin{displaymath}
  \begin{array}{rcccccccl}
    & \ldots
    & f^{-2\prime}(\theta)
    & f^{-1\prime}(\theta)
    & f(\theta)
    & f^{1\prime}(\theta)
    & f^{2\prime}(\theta)
    & \ldots
    &
    \\
    \mbox{integration}
    & \leftarrow
    & \updownarrow
    & \updownarrow
    & \updownarrow
    & \updownarrow
    & \updownarrow
    & \rightarrow\raisebox{-1.9ex}{\rule{0em}{5ex}}
    & \mbox{differentiation}
    \\
    & \ldots
    & w^{-2\ldot}(z)
    & w^{-1\ldot}(z)
    & w(z)
    & w^{1\ldot}(z)
    & w^{2\ldot}(z)
    & \ldots
    & .
  \end{array}
\end{displaymath}

\noindent
Note that there may be singular points of the inner analytic functions on
the unit circle, and that, as we discussed in~\cite{FTotCPII}, all the
inner analytic functions within a chain have exactly the same set of
singular points. Since any analytic function is both differentiable and
integrable, and its derivative and primitive are equally analytic, in
exactly the same domain as the original function, no point of singularity
ever appears or vanishes in processes of differentiation and integration.
Only the nature of the singular points, namely their degrees of softness
or degrees of hardness, changes due to the operations of logarithmic
differentiation and logarithmic integration.

In order to simplify the discussion that follows, let us establish some
limitations on the singularities that may be present. First, in the spirit
of the limitation on the inner analytic functions expressed by
Equation~(\ref{limitcond2}), let us eliminate from this discussion the
cases in which there are infinitely hard singularities, such as essential
singularities. Second, if the inner analytic functions in a chain have any
removable singularities, in what follows we will assume that they have
been removed. We will also treat separately any infinitely soft
singularities that may eventually exist in a chain. We are thus left only
with singularities at which some real function in the chain is
differentiable while its derivative is not. Finally, for simplicity of
argumentation, let us limit the following discussion to the case in which
there is a finite number of such singularities, which are then all
isolated singularities.

We may easily establish a general classification of all the
integral-differential chains, in a very simple way, by noting that if one
of the real objects in a chain is a $C^{\infty}$ zero-average real
function, then all the real objects in that chain must also be
$C^{\infty}$ zero-average real functions. On the other hand, if one of the
real objects is {\em not} a $C^{\infty}$ zero-average real function, then
none of them can be $C^{\infty}$. Therefore, either a chain contains only
$C^{\infty}$ zero-average real functions on its real sector, or it
contains none. We will call a chain that contains only $C^{\infty}$
zero-average real functions a $C^{\infty}$ chain. Note that the inner
analytic functions in any chain are all, of course, $C^{\infty}$ within
the open unit disk, in the complex sense, and that this statement and the
associated classification refers only to the real functions obtained on
the unit circle, and to whether or not they are $C^{\infty}$ on that
circle, in the real sense. As was discussed in~\cite{FTotCPII}, the case
of $C^{\infty}$ chains is that in which the inner analytic functions have
no singularities at all on the unit circle, or have only infinitely soft
singularities, such as the one exemplified in one of the appendices
of~\cite{FTotCPII}.

The other class of chains, which is the complement of the class of
$C^{\infty}$ chains, is the class of chains that include at least one
singularity of the inner analytic functions over the unit circle, which is
not an infinitely soft singularity. Let us recall that we assume that all
removable singularities have been eliminated, so that there are none of
those around to cloud the issue. Given the limitations on the types of
singularities that we have established, it follows that this is a
singularity at which, as we travel along the chain in the differentiation
direction, we eventually pass from a differentiable zero-average real
function to one that is not differentiable. From this fact one can draw
two conclusions, if we start from the last differentiable real function in
the chain. On the one hand, it follows that in the integration direction
of the chain all subsequent real objects are also differentiable
zero-average real functions. On the other hand, however, it also follows
that in the differentiation direction we immediately encounter a
zero-average real function that has a point of non-differentiability.
Immediately after that, we will encounter a zero-average real function
that has, at that singular point, either a step-type discontinuity or a
divergence to infinity.

Let us consider first the case in which there is a divergence to infinity.
If there is a singular point where the resulting real function diverges to
infinity, then upon further differentiation the function will become
non-integrable around that point. For example, if the divergence around
the point $\theta_{1}$ behaves as $1/(\theta-\theta_{1})^{p}$ with
$0<p<1$, which is integrable, then the first derivative behaves as
$1/(\theta-\theta_{1})^{p'}$ with $1<p'<2$, which is not integrable.
Therefore, we see that there are at least some non-integrable real
functions, in the sense that the integrals of their absolute values over
the unit circle diverge to infinity, which are still representable by
inner analytic functions. These are the locally non-integrable functions
discussed in Section~\ref{SEClocalnon}. In such cases the Fourier
coefficients of the real function cannot be calculated directly on the
unit circle in the usual way, but the functions can still be defined by
means of the corresponding inner analytic function. These are the cases to
which the definition of the extended Fourier coefficients can be applied.

Let us consider now the case in which there is a step-type discontinuity.
Since the differentiation of a step-type discontinuity results in a delta
``function'', the next step in the differentiation direction will produce,
therefore, a combination of real objects which includes a delta
``function'' at that point. After that, further steps in the
differentiation direction will produce combinations of real objects which
include the delta ``function'' and its successive higher-order
derivatives. Therefore, every chain included in this sub-class contains as
part of its real sector the delta ``function'' and {\em all} its multiple
derivatives. We may conclude, therefore, that the set of singular real
objects examined in Section~\ref{SECgenfunc}, which is a subset of the set
of all generalized functions, is an inevitable companion of the set of
zero-average integrable real functions on the unit circle.

Both the zero-average integrable real functions and the set of radically
singular objects represented by the delta ``function'' and its derivatives
of all orders are inevitably integrated as part of one and the same
structure. In addition to this, all locally non-integrable real functions
that are derivatives of any finite order of integrable real functions are
also part of that structure. In this way we see that the
integral-differential chains show that all the real objects studied in
this paper, either normal or generalized functions, either regular or
singular, are closely integrated into a single overall structure.

\section{Conclusions}

We have shown that the correspondence between real functions and inner
analytic functions established in earlier papers can be greatly extended.
With a set of inner analytic functions that is still limited by the
restriction expressed in Equations~(\ref{limitcond1})
and~(\ref{limitcond2}), adopted in order to avoid the most severe
singularities at the unit circle, it is possible to define a much larger
set of objects on the unit circle than was previously thought. This
includes not only all integrable real functions, regardless of any parity
properties, but at least some radically singular objects such as
generalized functions, in the general spirit of the Schwartz theory of
distributions, as well as at least a fairly large class of plainly
non-integrable real functions.

This development lead to the introduction of an extended definition of the
Taylor-Fourier coefficients of a real function, that can be used even for
the class of non-integrable real functions just mentioned. Coupled with
the technique of singularity factorization, which we introduced before,
these extended coefficients have an algorithmic value in terms of the
representation of these singular functions by trigonometric series. In a
somewhat surprising way, the concept of integral-differential chains,
introduced before for a completely different reason, found use here as a
way to systematize at least a large part of this large set of real
objects. In fact, it was instrumental for the very realization that it was
possible to extend the structure to non-integrable real functions, and
eventually led to the definition of the extended Taylor-Fourier
coefficients given in Equation~(\ref{extendefcoef}).

The set of inner analytic functions, as defined in the new standard way
proposed in this paper, has the properties of a vector space with real
scalars, in such a way that the addition of any two real objects on the
unit circle always corresponds to the addition of the two corresponding
inner analytic functions. In this way, arbitrary linear combinations of
any of the real objects discussed here also correspond to inner analytic
functions. This includes linear combinations mixing normal integrable
functions, locally non-integrable functions and singular generalized
functions. They can all be treated on the same footing. As a simple
example, generalized functions with several singular points can be
obtained by just adding the inner analytic functions corresponding to
several generalized functions with one singular point each.

Also in a somewhat surprising turn of events, an interesting connection
with the Dirichlet problem on the unit disk suggested itself in a rather
forceful way. This led at once to simple and easy extensions of the basic
theorems of existence and uniqueness of solutions to that problem, as well
as to the proposition of new types of boundary value problems, with what
we named punctured boundary conditions, which allow for the presence of
almost arbitrary singular points at the boundary. One finds that the
Dirichlet problem on the unit circle has a unique solution for any
integrable boundary condition, which can be discontinuous and even
unbounded at a finite set of singular points. The same is true for
singular boundary conditions consisting of generalized functions or
locally non-integrable functions.

The problem of the complete characterization of all the real objects on
the unit circle corresponding to the set of inner analytic functions
considered here remains open. One interesting question is whether or not
all Lebesgue-measurable but non-integrable real functions are included in
this structure, as limits of inner analytic functions to the unit circle.
Another is whether or not there are other radically singular objects, such
as but other than the Dirac delta ``function'' and its derivatives of all
orders, that can also be found as part of the structure.

\appendix

\section{Equivalence of Integrability and Local
  Integrability}\label{APPequivinteg}

Consider the set of Lebesgue-measurable real functions $f(\theta)$ defined
on $[-\pi,\pi]$. Let us recall that within this set the conditions of
integrability and of absolute integrability are equivalent conditions.
The condition of local integrability is defined as integrability on all
closed sub-intervals of $[-\pi,\pi]$. Since this includes $[-\pi,\pi]$
itself, it is immediate that the condition of local integrability implies
the condition of integrability.

Now consider an arbitrary closed sub-interval $[a,b]$ of $[-\pi,\pi]$.
Since the condition of integrability implies the condition of absolute
integrability, if $f(\theta)$ is integrable on $[-\pi,\pi]$, then it is
absolutely integrable there. But according to the properties of the
operation of integration we have that

\begin{displaymath}
  \int_{-\pi}^{\pi}d\theta\,
  |f(\theta)|
  =
  \int_{-\pi}^{a}d\theta\,
  |f(\theta)|
  +
  \int_{a}^{b}d\theta\,
  |f(\theta)|
  +
  \int_{b}^{\pi}d\theta\,
  |f(\theta)|.
\end{displaymath}

\noindent
If $f(\theta)$ is an absolutely integrable function, then in this relation
the number in the left-hand side is a finite positive real number, while
the right-hand side is the sum of three positive real numbers. It then
follows that these three real numbers must all be finite, and therefore
that $f(\theta)$ is absolutely integrable on $[a,b]$, since the integral
of its absolute value on that interval is a finite real number. Now, since
the condition of absolute integrability implies the condition of
integrability, it follows that $f(\theta)$ is integrable on $[a,b]$.

Therefore, since $f(\theta)$ is thus seen to be integrable on an arbitrary
closed sub-interval $[a,b]$ of $[-\pi,\pi]$, it follows that it is
integrable on {\em all} such closed sub-intervals, and is therefore
locally integrable. The final conclusion is that the conditions of
integrability and of local integrability are equivalent conditions in the
space of all Lebesgue-measurable real functions defined within
$[-\pi,\pi]$. Hence in this context all three conditions, that of local
integrability, that of absolute integrability and that of integrability,
are equivalent conditions.

\section{Calculation of some Inner Analytic Functions}\label{APPcalcoefs}

In this section we will derive a simple algebraic recursion relation for
all the inner analytic functions associated to the multiple derivatives of
the delta ``function''. Starting from the inner analytic function that
corresponds to the delta ``function'', given in
Equation~(\ref{innerdelta}) of the text,

\begin{displaymath}
  w(z,z_{1})
  =
  \frac{1}{2\pi}
  -
  \frac{1}{\pi}\,
  \frac{z}{z-z_{1}}, 
\end{displaymath}

\noindent
and from the differential recursion relation in terms of logarithmic
derivatives that relates each inner analytic function in the sequence to
the next one,

\begin{displaymath}
  w_{\delta^{(n+1)\prime}}(z,z_{1})
  =
  \ii z\,
  \frac{d}{dz}w_{\delta^{n\prime}}(z,z_{1}),
\end{displaymath}

\noindent
we may in principle derive all the inner analytic functions in the
sequence. Doing this one gets for the first few cases of the inner
analytic functions, including the one for the original delta ``function'',
for which we use the notation
$w_{\delta^{0\prime}}(z,z_{1})=w_{\delta}(z,z_{1})$, the list of results
that follows,

\noindent
\begin{eqnarray*}
  w_{\delta^{0\prime}}(z,z_{1})
  & = &
  \frac{1}{2\pi}
  -
  \frac{1}{\pi\ii^{0}}\,
  \frac{z}{z-z_{1}},
  \\
  w_{\delta^{1\prime}}(z,z_{1})
  & = &
  -\,
  \frac{1}{\pi\ii^{1}}\,
  \frac
  {zz_{1}}
  {(z-z_{1})^{2}},
  \\
  w_{\delta^{2\prime}}(z,z_{1})
  & = &
  -\,
  \frac{1}{\pi\ii^{2}}\,
  \frac
  {z(z+z_{1})z_{1}}
  {(z-z_{1})^{3}},
  \\
  w_{\delta^{3\prime}}(z,z_{1})
  & = &
  -\,
  \frac{1}{\pi\ii^{3}}\,
  \frac
  {z\left(z^{2}+4zz_{1}+z_{1}^{2}\right)z_{1}}
  {(z-z_{1})^{4}},
  \\
  w_{\delta^{4\prime}}(z,z_{1})
  & = &
  -\,
  \frac{1}{\pi\ii^{4}}\,
  \frac
  {z\left(z^{3}+11z^{2}z_{1}+11zz_{1}^{2}+z_{1}^{3}\right)z_{1}}
  {(z-z_{1})^{5}},
  \\
  w_{\delta^{5\prime}}(z,z_{1})
  & = &
  -\,
  \frac{1}{\pi\ii^{5}}\,
  \frac
  {z\left(z^{4}+26z^{3}z_{1}+66z^{2}z_{1}^{2}+26zz_{1}^{3}+z_{1}^{4}\right)z_{1}}
  {(z-z_{1})^{6}}.
\end{eqnarray*}

\noindent
Instead of working out the details of these first few cases, let us
construct inductively a simple algebraic recursion relation which gives
the general form of these inner analytic functions. We may write, in terms
of the variable $\chi=z/z_{1}$, for $n>0$, that

\begin{displaymath}
  w_{\delta^{n\prime}}(\chi)
  =
  -\,
  \frac{1}{\pi\ii^{n}}\,
  \frac{\chi P_{n-1}(\chi)}{(\chi-1)^{n+1}},
\end{displaymath}

\noindent
where $P_{n-1}(\chi)$ is a polynomial of order $n-1$ on $\chi$. The first
few polynomials, for $n>0$, are therefore given by

\noindent
\begin{eqnarray*}
  P_{0}(\chi)
  & = &
  1,
  \\
  P_{1}(\chi)
  & = &
  1+\chi,
  \\
  P_{2}(\chi)
  & = &
  1+4\chi+\chi^{2},
  \\
  P_{3}(\chi)
  & = &
  1+11\chi+11\chi^{2}+\chi^{3},
  \\
  P_{4}(\chi)
  & = &
  1+26\chi+66\chi^{2}+26\chi^{3}+\chi^{4}.
\end{eqnarray*}

\noindent
Since we have the differential recurrence relation for the inner analytic
function $w_{\delta^{n\prime}}(z,z_{1})$, in terms of logarithmic
derivatives,

\begin{displaymath}
  w_{\delta^{(n+1)\prime}}(\chi)
  =
  \ii\chi\,
  \frac{d}{d\chi}w_{\delta^{n\prime}}(\chi),
\end{displaymath}

\noindent
we may derive a related differential recurrence relation for the
polynomials,

\noindent
\begin{eqnarray*}
  -\,
  \frac{1}{\pi\ii^{n+1}}\,
  \frac{\chi P_{n}(\chi)}{(\chi-1)^{n+2}}
  & = &
  \ii\chi\,
  \frac{d}{d\chi}
  \left[
    -\,
    \frac{1}{\pi\ii^{n}}\,
    \frac{\chi P_{n-1}(\chi)}{(\chi-1)^{n+1}}
  \right]
  \\
  & = &
  -\,
  \frac{\ii\chi}{\pi\ii^{n}}
  \left[
    \frac{P_{n-1}(\chi)}{(\chi-1)^{n+1}}
    +
    \frac{\chi P'_{n-1}(\chi)}{(\chi-1)^{n+1}}
    -
    \frac{(n+1)\chi P_{n-1}(\chi)}{(\chi-1)^{n+2}}
  \right]
  \\
  & = &
  \frac{\chi}{\pi\ii^{n+1}}\,
  \frac{-(n\chi+1)P_{n-1}(\chi)+\chi(\chi-1)P'_{n-1}(\chi)}{(\chi-1)^{n+2}},
\end{eqnarray*}

\noindent
where $P'_{n}(\chi)$ is the derivative of $P_{n}(\chi)$ with respect to
$\chi$. We get therefore for the polynomials the differential recurrence
relation

\begin{displaymath}
  P_{n}(\chi)
  =
  (n\chi+1)P_{n-1}(\chi)-\chi(\chi-1)P'_{n-1}(\chi).
\end{displaymath}

\noindent
In order to solve this recurrence relation we write the polynomials
explicitly,

\noindent
\begin{eqnarray*}
  P_{n}(\chi)
  & = &
  \sum_{i=0}^{n}
  A_{n,i}\chi^{i}
  \;\;\;\Rightarrow
  \\
  P_{n-1}(\chi)
  & = &
  \sum_{i=0}^{n-1}
  A_{n-1,i}\chi^{i}
  \;\;\;\Rightarrow
  \\
  P'_{n-1}(\chi)
  & = &
  \sum_{i=0}^{n-1}
  iA_{n-1,i}\chi^{i-1},
\end{eqnarray*}

\noindent
so that we have for the differential recurrence relation

\noindent
\begin{eqnarray*}
  \sum_{i=0}^{n}
  A_{n,i}\chi^{i}
  & = &
  n\chi P_{n-1}(\chi)
  +
  P_{n-1}(\chi)
  -
  \chi^{2}P'_{n-1}(\chi)
  +
  \chi P'_{n-1}(\chi)
  \\
  & = &
  n
  \sum_{i=0}^{n-1}
  A_{n-1,i}\chi^{i+1}
  +
  \sum_{i=0}^{n-1}
  A_{n-1,i}\chi^{i}
  -
  \sum_{i=0}^{n-1}
  iA_{n-1,i}\chi^{i+1}
  +
  \sum_{i=0}^{n-1}
  iA_{n-1,i}\chi^{i}
  \\
  & = &
  n
  \sum_{i=1}^{n}
  A_{n-1,i-1}\chi^{i}
  -
  \sum_{i=1}^{n}
  (i-1)A_{n-1,i-1}\chi^{i}
  +
  \sum_{i=0}^{n-1}
  A_{n-1,i}\chi^{i}
  +
  \sum_{i=0}^{n-1}
  iA_{n-1,i}\chi^{i}.
\end{eqnarray*}

\noindent
If we now read out the special case $i=0$ we simply get
$A_{n,0}=A_{n-1,0}$, which is valid for all $n>0$. Since $P_{0}(\chi)=1$
implies that $A_{0,0}=1$, we have $A_{n,0}=1$ for all $n$. The special
case $i=n$ gives the equally simple result $A_{n,n}=A_{n-1,n-1}$, which is
also valid for all $n>0$. Since $A_{0,0}=1$, we have $A_{n,n}=1$ for all
$n$. For all the other values of $i$ we have the ordinary algebraic
recurrence relation

\begin{displaymath}
  \sum_{i=1}^{n-1}
  A_{n,i}\chi^{i}
  =
  \sum_{i=1}^{n-1}
  [(i+1)A_{n-1,i}+(n-i+1)A_{n-1,i-1}]\chi^{i},
\end{displaymath}

\noindent
so that we get the two-index algebraic recurrence relation for the
coefficients of the polynomials, given in Equation~(\ref{algebrecur}) of
the text,

\begin{displaymath}
  A_{n,i}
  =
  (i+1)A_{n-1,i}+(n-i+1)A_{n-1,i-1},
\end{displaymath}

\noindent
for $0<i<n$. Observe that this relation implies that all coefficients
$A_{n,i}$ will be integers. The special cases examined separately before
can now be merged back in if we agree that $A_{n,i}=0$ for $i<0$ and for
$i>n$. In the way of illustration, we may write some particular cases of
this recurrence relation, choosing particular values of $i$. For $i=1$ and
$i=n-1$ we get, respectively

\noindent
\begin{eqnarray*}
  A_{n,1}
  & = &
  nA_{n-1,0}
  +
  2A_{n-1,1},
  \\
  A_{n,n-1}
  & = &
  nA_{n-1,n-1}
  +
  2A_{n-1,n-2},
\end{eqnarray*}

\noindent
which are recurrence relations in terms of $n$ alone. Alternatively,
starting from the coefficient $A_{n,0}=1$, which has the value shown for
all $n$, and using the algebraic recurrence relation for a fixed value of
$n$, we may derive all the $n+1$ coefficients $A_{n,i}$ of the polynomial
$P_{n}(\chi)$. In this way we may construct the table of values of
$A_{n,i}$ shown in Table~\ref{thetable} of the text.

\end{document}